\documentclass[11pt]{article}
\usepackage[numbers,sort&compress]{natbib}
\usepackage{appendix}
\usepackage{enumerate}
\usepackage{amscd}
\usepackage{amsmath}
\usepackage{latexsym}
\usepackage{amsfonts}
\usepackage{amssymb}
\usepackage{amsthm}
\usepackage{verbatim}
\usepackage{mathrsfs}
\usepackage{enumerate}
\usepackage{hyperref}

\theoremstyle{plain}\newtheorem{definition}{Definition}[section]
\theoremstyle{definition}\newtheorem{theorem}{Theorem}[section]
\theoremstyle{plain}\newtheorem{lemma}[theorem]{Lemma}
\theoremstyle{plain}\newtheorem{coro}[theorem]{Corollary}
\theoremstyle{plain}
\theoremstyle{remark}\newtheorem{remark}{Remark}[section]
\usepackage{xcolor}
\usepackage[normalem]{ulem}

\newcommand{\norm}[1]{\left\|#1\right\|}

\newcommand{\B}{\Big}

\newcommand{\s}{\mathrm{div}}

\newcommand{\be}{\begin{equation}}
\newcommand{\ee}{\end{equation}}
 \newcommand{\ba}{\begin{aligned}}
 \newcommand{\ea}{\end{aligned}}

  \newcommand{\f}{\frac}
    
  \newcommand{\ben}{\begin{enumerate}}
   \newcommand{\een}{\end{enumerate}}

\newcommand{\Rmnum}[1]{\expandafter\@slowromancap\romannumeral #1@}

\allowdisplaybreaks

\numberwithin{equation}{section}
\begin{document}
\title{Helicity conservation for compressible barotropic Euler flows: regularity criteria and a momentum-based approach}
\author{ Pu Wang\footnote{School of Mathematics and Statistics,
		Henan University,
		Kaifeng, 475004,
		P. R. China. Email: wangpu.math@henu.edu.cn, wangpu980809@163.com }, \,~
Xiang Xia\footnote{ School of Mathematics and Statistics,
      	Henan University,
      	Kaifeng, 475004,
      	P. R. China. Email: xiaxiang0379@163.com},
      \, ~  and  ~  Yulin Ye\footnote{Corresponding author.   School of Mathematics and Statistics,
      	Henan University,
      	Kaifeng, 475004,
      	P. R. China. Email: ylye@vip.henu.edu.cn} }
\date{}
\maketitle
\begin{abstract}
In this paper, we study helicity conservation for weak solutions to the compressible barotropic Euler equations. Using generalized Constantin-E-Titi-type commutator estimates that fully exploit the  kinetic energy, we derive new sufficient conditions for helicity conservation across a broader range of regularity exponents. These conditions reveal the interplay between temporal, spatial integrability and differentiability of the velocity and density, as well as the influence of velocity-density coupling on helicity conservation. A key feature of our approach is that helicity conservation is established directly from the momentum equation, without invoking the vorticity equation.

\end{abstract}

\noindent {\bf MSC(2020):}\quad 76N10, 35Q31, 35L65\\\noindent
{\bf Keywords:} compressible Euler equations;  weak solutions, helicity conservation;   
\section{Introduction}
In the present paper, we consider the  following  three-dimensional ideal compressible barotropic fluid, the model of which can be described as follows:
\begin{equation}\label{CE}
\begin{cases}
\partial_t \rho + \s(\rho v) = 0,\\
\partial_t(\rho v) + \s(\rho v\otimes v) + \nabla \pi(\rho) = 0,\\
\rho|_{t=0}=\rho_0,~\rho v|_{t=0}=\rho_0 v_0,
\end{cases}
\qquad (t,x)\in (0,T)\times\Omega,
\end{equation}
where  the scalar function $\rho$ and the unknown vector $v $ represent   the density and  velocity of the flow,  respectively. $\pi(\rho)$ stands for the   pressure of the fluid. Here, a barotropic fluid refers to the case where pressure depends solely on density. For the sake of simplicity, we consider the problem in a bounded domain with periodic boundary condition, namely, $\Omega=\mathbb{T}^3$.

Formally, when $\rho=constant$, then the compressible Euler equations \eqref{CE} reduce to the ideal incompressible Euler equations
  \begin{equation}\left\{\begin{aligned}\label{Euler}
		&v_{t}+v\cdot \nabla v+\nabla \pi=0,\\
		&\text{div}\, v=0,\\
		&v|_{t=0}=v_{0}(x),
	\end{aligned}\right.\end{equation}
which is widely used to model the turbulence when the Reynolds number is large enough. As is known to all, there exist two physically quadratic conserved  quantities for the regular solutions of the Euler system \eqref{Euler}:
\begin{itemize}
	\item  Energy  conservation
	\be\label{ec}
	\f12\int_{\Omega}|v(t,x)|^{2}dx=\f12\int_{\Omega}|v(0,x)|^{2}dx;
	\ee
	\item Helicity preservation
	\be\label{hc}
	\f12\int_{\Omega} \omega(t,x)\cdot v(t,x)  dx=\f12\int_{\Omega} \omega(0,x)\cdot v(0,x)  dx,
	\ee
\end{itemize}
where $\omega=curl v$ is the vorticity of the fluid.
The energy of fluid is
the
core  in  Kolmogorov's 1941 turbulence theory   and the Onsager conjecture. In particular, 
Onsager in \cite{[Onsager]} asserted   that weak solutions of the ideal incompressible Euler equations with
H\"older continuity exponent $\alpha>\frac{1}{3}$ do conserve energy and that turbulent
or anomalous dissipation occurs when $\alpha\leq \frac{1}{3}$ in 1949.  For the progresses  in both two directions of the Onsager's conjecture, the readers can be referred to  \cite{[B],[BG1],[GKN],[LS],[BV],[GR],[Eyink],[CET],[DR],[CCFS],[Shvydkoy2009],[FW2018],[DS0],[Isett],[ADSW],[RDII],[B]}. Besides, 
Onsager'  conjecture involving energy conservation
for the compressible Euler   equations \eqref{CE} with vacuum can be found in \cite{[ADSW],[BGSTW],[EGSW]}.

Actually, the description of turbulence not only rests on the cascades of energy, but also the helicity. The interaction of the transfer between energy and helicity plays a critical  role in determining the direction of their  cascade
(see \cite{[HLY]}). Historically, the concept of the
helicity  in an inviscid compressible barotropic fluid was proposed  by Moffatt in \cite{[Moffatt]}, related discussions can also be found in \cite{[Moreau]}. As pointed in \cite{[Moffatt],[MT]},
helicity is important at a fundamental level in relation to flow kinematics
because it admits topological interpretation in relation to the linkage or
linkages of vortex lines of the flow.
 Consequently, considerable progress has been made in understanding helicity conservation for weak solutions to the incompressible Euler equations. Firstly,  Chae \cite{[Chae],[Chae1]} derived sufficient conditions for helicity conservation involving the vorticity $\omega$. More precisely, helicity conservation holds if $\omega\in L^{3}(0,T;B^{\alpha}_{\f95,\infty})$ with $\alpha>\f13$, or if $v\in L^{r_{1}}(0,T;B^{\alpha}_{\f92,q})$ and $\omega\in L^{r_{2}}(0,T;B^{\alpha}_{\f95,q})$, where $\alpha>\f13$, $q\in[2,\infty]$, $r_1\in[2,\infty]$, $r_2\in[1,\infty]$, and $\frac{2}{r_1}+\frac{1}{r_2}=1$.  In particular, Cheskidov, Constantin, Friedlander and Shvydkoy in \cite{[CCFS]} established the helicity conservation class based on the velocity $v$ in  Onsager's critical space
 $   L^{3}(0,T;B^{\f23}_{3,c(\mathbb{N})})$. From a mathematical viewpoint, since the nonlinear term $\int_{\Omega} (v\cdot \nabla v)\cdot \omega dx$, which arises in the helicity computation for the Euler equations \eqref{Euler}, contains two derivatives and three copies of the velocity $v$, it seems that
 $2/3$ is the natural  Onsager exponent for helicity conservation. 
In the framework of Sobolev spaces,  De Rosa \cite{[De Rosa]}  further established helicity conservation under the assumptions 
$v\in L^{2r}(0,T;W^{\theta, 2p})$ and  $ \omega\in L^{\kappa}(0,T;W^{\alpha, q})$ with  $\f{1}{p}+\f1q=\f{1}{r}+\f{1}{\kappa}$  and $2\theta+\alpha\geq1$.  Besides, Boutros and Titi \cite{[BT]} also characterized the helicity defect arising from the lack of smoothness of weak solutions. Inversi and Sorella \cite{[IS]} investigated the global helicity balance for the Euler equations \eqref{Euler} in bounded domains, accounting for boundary contributions involving the vorticity, velocity, and pressure under suitable regularity assumptions. It is worth emphasizing that the aforementioned results, together with those in \cite{[LZ]}, are confined to the incompressible setting.

However,   there   has been very few literature  concerning the   helicity conservation  of weak solutions for
the    compressible case. A key distinction between incompressible and compressible Euler systems is that one cannot directly derive a vorticity equation from the momentum equation in the compressible setting, which is  a step central to analyzing helicity conservation in the incompressible case. This difficulty arises from the low regularity of \((\rho,v)\), the strong coupling between density \(\rho\) and velocity $v$, and the possible presence of vacuum states.  Recently, the authors in \cite{[Ye2]} observed that for any smooth solutions $(\rho, v)$ of the compressible Euler equations \eqref{CE} with $0<c_1\leq \rho\leq c_2<\infty$, dividing the both sides of the momentum equation $\eqref{CE}_2$ by $\rho$, then the compressible Euler equations can be  reformulated as
\be\left\{\ba\label{rCEuler1}
&\rho_t+\s (\rho v)=0, \\
& v_{t} + v\cdot\nabla v+\nabla
P(\rho )=0,
\ea\right.\ee
where  $P(\rho )= \int_{1}^{\rho}\f{\pi'(s)}{s}ds$.
Then the momentum equation  in \eqref{rCEuler1} is the same as that for
the incompressible case. Based on this observation, the authors in \cite{[Ye2]} considered this issue in two steps and proved that if the weak solutions satisfy
	\begin{equation}\label{a1}\ba
	&0<c_1\leq \rho \leq c_2<\infty,~	~\rho \in L^3(0,T;\dot{B}^{\frac{1}{3}}_{3,c(\mathbb{N})}),~\rho v \in L^3(0,T;\dot{B}^{\frac{1}{3}}_{3,\infty}),\\
	&v \in L^3(0,T;B^{\frac{1}{3}}_{3,c(\mathbb{N})}),\
	\pi\in C^2[c_1,c_2],
	\ea\end{equation}
	then $\omega\in L^{3}(0,T;\dot{B}^{\frac{1}{3}}_{3,\infty})$
	or $\omega \in L^3(0,T;L^3(\mathbb{T}^3))$ ensure the  helicity identity for the compressible Euler equations \eqref{CE}. It should be noted in \cite{[Ye2]} the velocity is required to belong to $v\in L^3(0,T;\dot{B}^{\f43}_{3,c(\mathbb{N})})$ or $v\in L^3(0,T;{B}^{1}_{3,3})$, which has a gap compared with the critical space $v\in L^3(0,T;{B}^{\f23}_{3,c(\mathbb{N})})$. The main reasons are twofold. First,  the proof in \cite{[Ye2]} is divided into two steps and it mainly relies on the research strategy of the incompressible Euler equations \eqref{Euler}. Second one is that the vector identity $v\cdot \nabla v= \s (v\otimes v)$ is nolonger valid in the case of the compressible Euler equations \eqref{CE} due to the absence of the divergence-free condition. This naturally raises the question of whether helicity conservation for the compressible Euler equations \eqref{CE} can be established under the Onsager-critical regularity assumption $v\in L^3(0,T;B^{\f23}_{3,c(\mathbb{N})})$, as in the incompressible setting. Another important question is how density variations and their coupling with the velocity field affect helicity conservation in compressible flows. Motivated by these questions, the present paper aims to establish refined regularity criteria for helicity conservation while elucidating the role of density in the compressible setting.  Now, we formulate our main result as follows.
\begin{theorem}\label{the 1.1}
Let $(\rho,v)$ be a weak solution to the three-dimensional compressible barotropic  Euler system \eqref{CE} in the Definition   \ref{ceulerdefi}. Suppose $\omega\in C([0,T];L^{\f{3}{2}}(\mathbb{T}^3))$, $v\in C([0,T];L^{3}(\mathbb{T}^3))$. 
Assume further that $0<c_1\leq \rho \leq c_2<\infty$ and that one of the following conditions holds:
\begin{enumerate}[(i)]
\item for $\frac{2}{3}<\alpha<1$, $\rho v, v\in L^{\frac{2}{\alpha}}\bigl(0,T;B^{\alpha}_{\frac{4}{2-\alpha},c(\mathbb{N})}(\mathbb{T}^3)\bigr)$,\newline $\rho \in L^{\frac{7}{4(1-\alpha)}}\bigl(0,T; B^{\max\{2-2\alpha, \frac{2-\alpha}{3}\}}_{\infty,\infty}(\mathbb{T}^3)\bigr)$ and $\pi\in C^2(0,\infty)$;
\item $\rho v, v\in L^{3}\bigl(0,T;B^{\frac{2}{3}}_{3,c(\mathbb{N})}(\mathbb{T}^3)\bigr)$, $\rho \in L^{\f{21}{4}}\bigl(0,T; B^{\frac{2}{3}}_{\infty,\infty}(\mathbb{T}^3)\bigr)$ and $\pi\in C^1(0,\infty)$;
\item $v\in L^{2}\bigl(0,T;W^{1,4}(\mathbb{T}^3)\bigr)$, $\rho \in L^{\infty}\bigl(0,T; B^{\frac{1}{3}}_{\infty,c(\mathbb{N})}(\mathbb{T}^3)\bigr)$ and $\pi\in C^2(0,\infty)$;
\item $v\in L^{2}\bigl(0,T;W^{1,4}(\mathbb{T}^3)\bigr)$, $\rho v\in L^{1}\bigl(0,T;W^{1,1}(\mathbb{T}^3)\bigr)$, $\rho\in L^{\infty}\bigl(0,T;B^{\frac{1}{3}}_{\infty,\infty}(\mathbb{T}^3)\bigr)$ and $\pi\in C^2(0,\infty)$.
\end{enumerate}
Then the helicity of  the weak solutions to the compressible Euler equations is globally conserved, that is, for all $0<t<T$,
$$	\int_{\mathbb{T}^3} \omega(t,x)\cdot v(t,x)  dx=\int_{\mathbb{T}^3} \omega(0,x)\cdot v(0,x)  dx.$$

\end{theorem}
\begin{remark}\label{rem1}
We establish a sufficient condition for helicity conservation of the compressible Euler equations under a broader range of exponent scales. This result not only clarifies how the time integrability, spatial integrability, and differential regularity of the velocity and density fields relate to one another, but also illuminates how the coupling between velocity and density contributes to helicity conservation for weak solutions. To the authors’ best knowledge, this appears to be the first result revealing the interplay between density and velocity in the helicity conservation law for the compressible Euler equations.
\end{remark}
\begin{remark}\label{rem4}
The velocity field $v$ requires  the regularity $L^3\bigl(0,T; B^{\f{2}{3}}_{3,c(\mathbb N)}(\mathbb{T}^3)\bigr)$ to guarantee the helicity conservation in the compressible Euler equations \eqref{CE}, 
which is consistent with the helicity conservation criterion in the case of the incompressible Euler equations \eqref{Euler} in \cite{[CCFS]}.
Roughly speaking, we extend the critical sufficient condition for the helicity conservation from the incompressible Euler equations to the compressible case. Moreover, we give an affirmative answer to the open problem proposed in \cite{[Ye2]}.
\end{remark}

\begin{remark}\label{rem3}
We remark that, in the endpoint case $\alpha=\frac{2}{3}$ of Theorem~\ref{the 1.1} (ii), the regularity requirement on the pressure is weakened compared with that in the subcritical case (i). Similarly, in the endpoint case $\alpha=1$ of Theorem~\ref{the 1.1} (iv), the regularity assumption on the momentum $\rho v$ is further weakened by splitting the domain into the regions where $v=0$ and $v\neq 0$ and treating them separately.
\end{remark}
\begin{remark}
	In the endpoint case $\alpha=1$ of Theorem~\ref{the 1.1} (iii), we eliminate the additional regularity assumption on the momentum $\rho v$ at the expense of imposing a slightly stronger regularity condition on the density, namely, $\rho\in L^{\infty}\bigl(0,T;B^{\frac{1}{3}}_{\infty,c(\mathbb{N})}(\mathbb{T}^3)\bigr)$.
\end{remark}
\begin{remark}\label{rem6}
Compared with all previous results on helicity conservation, which rely heavily on the vorticity equation, we develop a new framework to analyze helicity using only the momentum equation.	
\end{remark}
\begin{remark}\label{rem2}
	It is worth remarking  that the assumptions $\omega\in C([0,T];L^{\f{3}{2}}(\mathbb{T}^3))$ and  $v\in C([0,T];L^{3}(\mathbb{T}^3))$ are employed to make sure the helicity is well-defined. However, these assumptions are not necessary and, naturally, not unique.	
\end{remark}
We give some comments on the proof of Theorem \ref{the 1.1}. Owing to the rate of change of the helicity
\be\label{rateofhelicity}\f{d}{dt}\int_{\mathbb{T}^3} v^\varepsilon\cdot(\nabla \times v)^\varepsilon\, dx= \int_{\mathbb{T}^3} v_t^\varepsilon \cdot (\nabla \times v)^\varepsilon \,dx+\int_{\mathbb{T}^3} v^\varepsilon \cdot (\nabla \times v)_t^\varepsilon \,dx,
\ee
it is natural to apply both the velocity equations and the vorticity equations in the following computation (see \cite{[Chae],[Chae1],[De Rosa],[Ye2],[Moffatt],[MT],[BT],[IS],[Ye]}).
For the  incompressible case, the vorticity equations can be obtained by taking the $curl$ operator on both sides of the momentum equation.
 However,  this strategy breaks down in the compressible setting, due to the strong coupling between the density and the velocity in the momentum equation. To overcome this difficulty, with the help of the vector identity
 $$ {E}\cdot(\nabla\times  {F})=\text{div}( {F}\times {E})+ {F}\cdot(\nabla\times {E}),$$ we observe that the rate of change of the helicity \eqref{rateofhelicity} can be reformulated as
 \be\label{eqrateofhelicity}
 \f{d}{dt}\int_{\mathbb{T}^3} v^\varepsilon\cdot(\nabla \times v)^\varepsilon\, dx= 2\int_{\mathbb{T}^3} v_t^\varepsilon \cdot (\nabla \times v)^\varepsilon \,dx.
\ee
 The advantage of equation \eqref{eqrateofhelicity} is that the vorticity equations are not required straightforwardly.
Roughly speaking, we could choose
 $(\nabla \times v)^\varepsilon$ as the test function and compute the helicity of the compressible fluids in one step rather than two steps in \cite{[Ye2]}. It seems that the form of the rate of change of the helicity \eqref{rateofhelicity} is novel. After we finished this manuscript, we try our best to check previous works and find that similar identity  for  the magnetic helicity has been noticed in \cite{[Woltjer],[BV]}.

Notice that most aforementioned results are the rigid part of Onsager's conjecture concerning the helicity in both the incompressible and compressible flows.
A  breakthrough of
  flexible part of Onsager's conjecture for constructing the  weak
solutions with non-constant helicity was recently made by  Giri,  Kwon  and   Novack in \cite{[GKN2]}, where they
 constructed weak solutions with almost Onsager-critical regularity in $L^{3}$
  for given
generalized helicity and kinetic energy in the incompressible
ideal fluids \eqref{Euler} via extending the classical definition of the  helicity. Finally, we also would like to mention that Boutros  and  Gibbon addressed the removal of the barotropic condition
in helicity studies of the compressible Euler equations in \cite{[BG]}.

The rest of this paper is organized as follows.
In Section 2, we give  notations and some useful estimates, which play an important role in the follow-up study.
Section 3 is devoted to  the proof of helicity conservation of weak solutions to the 3-D compressible Euler equations. 
  
\section{Notations and some auxiliary lemmas} \label{section2}

 First, we introduce some notations used in this paper. We let $\varrho:\mathbb{R}^{3}\rightarrow \mathbb{R}$ be a standard mollifier, i.e., $\varrho(x)=C_0e^{-\frac{1}{1-|x|^2}}$ for $|x|<1$ and $\varrho(x)=0$ for $|x|\geq 1$, where $C_0$ is a constant such that $\int_{\mathbb{R}^3}\varrho (x) dx=1$. For $\varepsilon>0$, we define the rescaled mollifier $\varrho_\varepsilon(x)=\frac{1}{\varepsilon^3}\varrho(\frac{x}{\varepsilon})$ and for  any function $f\in L^1_{loc}(\mathbb{T}^3)$, its mollified version is defined as
\be\label{mollified}f^\varepsilon(x)=(f*\varrho_{\varepsilon})(x)=\int_{B_\varepsilon(0)}f(x-y)\varrho_\varepsilon(y)dy,\ \ x\in \mathbb{T}^3.\ee
Next, we introduce some function spaces.\\
$\mathbf{Sobolev \ spaces.}$  For $p\in [1,\,\infty]$, the notation $L^{p}(0,\,T;X)$ stands for the set of measurable functions on the interval $(0,\,T)$ with values in $X$ and $\|f(t,\cdot)\|_{X}$ belonging to $L^{p}(0,\,T)$. The classical Sobolev space $W^{k,p}(\mathbb{T}^3)$ is equipped with the norm $\|f\|_{W^{k,p}(\mathbb{T}^3)}=\sum\limits_{|\alpha| =0}^{k}\|D^{\alpha}f\|_{L^{p}(\mathbb{T}^3)}$.\\
{\bf Besov spaces.} For any  $0<s<1$ and $1\leq q \leq \infty$, we define the Besov semi-norm $\norm{f}_{\dot{B}^s_{q,\infty}(\mathbb{T}^3)}$ and Besov norm $\norm{f}_{B^s_{q,\infty}(\mathbb{T}^3)}$ of $f\in \mathcal{S}^{'}$ as
\begin{equation}\label{2.2}\|f\|_{\dot{B}_{q, \infty}^s(\mathbb{T}^3)}=\sup_{|y|>0}\f{\| f(x-y)-f(x)\|_{L^q(\mathbb{T}^3\cap (\mathbb{T}^3 +\{y\}))}}{|y|^s},
	\norm{f}_{B^s_{q,\infty}(\mathbb{T}^3)}=\norm{f}_{{L^q}(\mathbb{T}^3)}+\norm{f}_{\dot{B}^s_{q,\infty}(\mathbb{T}^3)},\end{equation}
where $\mathbb{T}^3 +\{y\}:=\{x\in \mathbb{T}^3| x=z+y, z\in \mathbb{T}^3\}$ with any $y\neq 0$ denoting the translation of the domain $\mathbb{T}^3 $ according to $y$. Furthermore, motivated by \cite{[CCFS],[Ye2]}, 
we denote ${B}^\alpha _{q,c(\mathbb{N})}$ with $\alpha\in (0,1)$ and $1\leq q \leq \infty$ as the class of all tempered distributions $f$ satisfying
\begin{equation}\label{2.3}
	\norm{f}_{{B}^\alpha _{q,\infty}(\mathbb{T}^3)}<\infty~ \text{and}~ 	\lim_{|y|\rightarrow 0} \f{\| f(x-y)-f(x)\|_{L^q(\mathbb{T}^3\cap (\mathbb{T}^3 +\{y\}))}}{|y|^\alpha}=0.
\end{equation}

Note that, combining the definitions of spaces ${B}^\alpha_{q,r}$  and $B^\alpha_{q,c(\mathbb{N})}$ above, we can deduce that
\begin{equation}\label{2.4}
{B}^\alpha_{q,r}(\mathbb{T}^3)	\hookrightarrow {B}^\alpha_{q,c(\mathbb{N})}(\mathbb{T}^3)\hookrightarrow B^\alpha_{q,\infty}(\mathbb{T}^3),
\end{equation} for any $1\leq q\leq \infty$, $1\leq r<\infty$ and $0<\alpha<1$.
Before going any further, for simplicity, we write
$$\|f\|_{L^p(0,T;L^q(\mathbb{T}^3))}=\|f\|_{L^p_t(L^q_x)}~~\text{and}~~\|f\|_{L^p(0,T;B^{\alpha}_{q,r}(\mathbb{T}^3))}=\|f\|_{L^p_t(B^\alpha_{q,r})}.$$

Next,  we present some  useful  lemmas concerning mollifier  and functions in Besov spaces, which play an important role in the proof.

\begin{lemma}\label{lem2.1}\cite{[Ye2],[Chae],[Chae1],[CET],[BGSTW]}
Let $\alpha, \beta\in (0,1)$,  $ p,q\in [1,\infty]$,  and $k\in \mathbb{N}^+$. Assume that  $f\in L^p(0,T;\dot{B}^\alpha_{q,\infty})$and  $g\in L^p(0,T;\dot{B}^\beta_{q,c(\mathbb{N})})$,  then as $\varepsilon\to 0$, there hold that
 \begin{enumerate}[(1)]
 \item $ \|f^{\varepsilon} -f \|_{L^{p}(0,T;L^{q})}\leq \text{O}(\varepsilon^{\alpha})$;
   \item   $ \|\nabla^{k}f^{\varepsilon}  \|_{L^{p}(0,T;L^{q})}\leq \text{O}(\varepsilon^{\alpha-k})$;
       \item $ \|g^{\varepsilon} -g \|_{L^{p}(0,T;L^{q})}\leq \text{o}(\varepsilon^{\beta})$;
   \item   $ \|\nabla^{k}g^{\varepsilon}  \|_{L^{p}(0,T;L^{q})}\leq \text{o}(\varepsilon^{\beta-k}),$
 \end{enumerate}
 where $f^{\varepsilon}$ and $g^{\varepsilon}$ are defined by \eqref{mollified}.
\end{lemma}
As a corollary of Lemma \ref{lem2.1}, we have the following properties.
\begin{coro}\label{coro}(\cite{[Shvydkoy2009],[Yu],[Ye2]})
	Let  $0<\alpha, \beta<1$ and  $1\leq q_1,q_2\leq \infty$. Assume that $f\in  \dot{B}^{\alpha}_{q_{1},\infty} $ and $g\in \dot{B}^{\beta}_{q_{2},c(\mathbb{N})} $, then there hold
	\begin{equation}\label{b5-2}\begin{aligned}
			&\|f(\cdot-y)-f(\cdot)\|_{L^{q_{1}}} \leq C|y|^{\alpha}\norm{f}_{\dot{B}^\alpha_{q_1,\infty}},\\
			&\|g(\cdot-y)-g(\cdot)\|_{L^{q_{2}}} \leq C|y|^{\beta}\norm{g}_{\dot{B}^\beta_{q_2,c(\mathbb{N})}}.
	\end{aligned}\end{equation}
	Moreover, as $|y|\to 0$, we have
	\begin{equation}\label{b5-22}\begin{aligned}
			&\|f(\cdot-y)-f(\cdot)\|_{L^{q_{1}}} = \text{O}(|y|^\alpha),\\
			&\|g(\cdot-y)-g(\cdot)\|_{L^{q_{2}}}   = \text{o}(|y|^\beta).
	\end{aligned}\end{equation}
\end{coro}
\begin{lemma} (\cite{[YeLi]})\label{lem2.3}
	Let $0<\alpha <1$, $0<\theta,\gamma\leq1$, $1\leq r,s,p_{i},q_i\leq\infty$ such that $\frac{1}{r}=\sum_{i=1}^4\frac{1}{p_i}$ and $\frac{1}{s}=\sum_{i=1}^4\frac{1}{q_i}$, for $i=1,2,3,4$. Suppose that  $f\in   L^{\theta p_{1}}(0,T;{B}_{\theta q_{1},c(\mathbb{N})}^{\alpha}(\mathbb{T}^3))\cap L^{(1-\theta)p_{2}}(0,T;L^{(1-\theta)q_{2}}(\mathbb{T}^3))$ and   $g\in  L^{\gamma p_{3}}(0,T;{B}_{\gamma q_{3},\infty}^{\beta}(\mathbb{T}^3))\cap  L^{(1-\gamma)p_{4}}(0,T;L^{(1-\gamma)q_{4}}(\mathbb{T}^3))$, then as $\varepsilon\rightarrow0,$ there holds
	\begin{align} \label{cet}
		\|(fg)^{\varepsilon}- f^{\varepsilon}g^{\varepsilon}\|_{L^r(0,T;L^s(\mathbb{T}^3))} \leq \text{o}(\varepsilon^{\theta\alpha+\gamma \beta }).	
	\end{align}

\end{lemma}
\begin{remark}
	Using the embedding relationship \eqref{2.4},  if we take $f\in   L^{\theta p_{1}}(0,T;B_{\theta q_{1},\infty}^{\alpha}(\mathbb{T}^3))$, then the result becomes $	\|(fg)^{\varepsilon}- f^{\varepsilon}g^{\varepsilon}\|_{L^r(0,T;L^s(\mathbb{T}^3))} \leq \text{O}(\varepsilon^{\theta\alpha+\gamma \beta }).$
\end{remark}
\begin{remark}
	When $\theta=\gamma =1$,  the result reduces to the classical Constantin-E-Titi type commutator estimates which can be found in \cite{[Ye2]}.
\end{remark}

 Furthermore, there also holds the following Sobolev-Besov mixed type commutator estimates.
\begin{lemma}\label{lem2.5}
Let $0<\alpha <1$ and $1\leq p,q< \infty$. Assume $f \in L^\infty(0,T;L^\infty(\mathbb{T}^3))$ and
$g \in L^p(0,T;\dot{B}^{\alpha}_{q,c(\mathbb N)}(\mathbb{T}^3))$, then there holds that
\be\label{2.8}
\|(fg)^{\varepsilon}- f^{\varepsilon}g^{\varepsilon}\|_{L^p(0,T;L^q(\mathbb{T}^3))} \leq  \text{o}(\varepsilon^{\alpha}),	\ \text{as}\ \varepsilon\to 0.
\ee
\end{lemma}
\begin{proof}
Thanks to the Constantin-E-Titi type identity observed in \cite{[CET]}, we have
\be\label{2.9}\ba
(fg)^\varepsilon-f^\varepsilon g^\varepsilon=R^\varepsilon-(f^\varepsilon-f)(g^\varepsilon-g),
\ea\ee
where
\begin{align*}
	R^\varepsilon(x,t):=\int_{\mathbb{R}^3}\varrho_\varepsilon(y)\B(f\left(t,x-y\right)-f(t,x)\B) \B(g(t,x-y)-g(t,x)\B)dy.
\end{align*}
On the one hand, by virtue of the H\"older inequality and   the
Minkowski inequality, we arrive at
\be\ba
\|R^\varepsilon\|_{L^q_x(\mathbb{T}^3)}
&\le
\int_{B_\varepsilon(0)}\varrho_\varepsilon(y)\,
\|f(t,\cdot-y)-f(t,\cdot)\|_{L^\infty}\,
\|g(t,\cdot-y)-g(t,\cdot)\|_{L^q_x}\,dy \notag\\
&\le
C\|f\|_{L^\infty}
\int_{B_\varepsilon(0)}\varrho_\varepsilon(y)\,
\f{\|g(t,\cdot-y)-g(t,\cdot)\|_{L^q_x}}{|y|^\alpha} |y|^\alpha \,dy\\
&\leq
C\varepsilon^\alpha\|f\|_{L^\infty}
 \|g\|_{\dot{B}^{\alpha}_{q, c(\mathbb{N})}}\int_{B_\varepsilon(0)}\varrho_\varepsilon(y)\,
 \,dy\\
 &\leq
 C\varepsilon^\alpha\|f\|_{L^\infty}\|g\|_{\dot{B}^{\alpha}_{q, c(\mathbb{N})}},
\ea\ee
which together with the H\"older inequality implies that
\be\label{2.10}\ba
\|R^\varepsilon\|_{L^p(0,T;L^q(\mathbb{T}^3))}&\leq
C\varepsilon^\alpha\|f\|_{L^\infty_t(L^\infty_x)}\|g\|_{L^p_t(\dot{B}^{\alpha}_{q, c(\mathbb{N})})}.
\ea\ee
On the other hand, by means of H\"older inequality and Lemma \ref{lem2.1}, we know that
\be\label{2.11}
\ba
\|(f^\varepsilon-f)(g^\varepsilon-g)\|_{L^q(\mathbb{T}^3)}&\leq \|f^\varepsilon-f\|_{L^\infty(\mathbb{T}^3)}\|g^\varepsilon-g\|_{L^q(\mathbb{T}^3)}\\
&\leq C\varepsilon^\alpha\|f\|_{L^\infty(\mathbb{T}^3)}\|g\|_{\dot{B}^{\alpha}_{q, c(\mathbb{N})}(\mathbb{T}^3)},
\ea\ee
which implies
\be\label{2.12}\ba
\|(f^\varepsilon-f)(g^\varepsilon-g)\|_{L^p(0,T;L^q(\mathbb{T}^3))}&\leq
C\varepsilon^\alpha\|f\|_{L^\infty_t(L^\infty_x)}\|g\|_{L^p_t(\dot{B}^{\alpha}_{q, c(\mathbb{N})})}.
\ea\ee
Then combining \eqref{2.9}, \eqref{2.10},\eqref{2.12} and the definition of Besov spaces $\dot{B}^\alpha_{q,c(\mathbb{N})}$, we get the desired estimate \eqref{2.8}.
\end{proof}

Next, we present the following estimate in the Besov-space setting, which plays a crucial role in the proof of the main result. 
\begin{lemma}\label{lem2.4}
For any $1\le q<\infty$ and $0<\alpha<1$, assume that
$0<c_1\le \rho \le c_2<\infty$ and $v\in B^\alpha_{q,c(\mathbb{N})}(\mathbb{T}^3)$,
 then there exists a constant $C>0$  independent of $\varepsilon$, such that
\begin{equation}\label{2.13}
\left\|
\nabla\left(\frac{(\rho v)^\varepsilon}{\rho^\varepsilon}\right)
\right\|_{L^q(\mathbb{T}^3)}
\le
C\,\varepsilon^{\alpha-1}\,
\|v\|_{B^\alpha_{q,c(\mathbb{N})}(\mathbb{T}^3)} .
\end{equation}
Furthermore, we also have
\be\label{2.14}\ba
\left\|
\nabla^2\left(\frac{(\rho v)^\varepsilon}{\rho^\varepsilon}\right)
\right\|_{L^q(\mathbb{T}^3)}
\le
C\,\varepsilon^{\alpha-2}\,
\|v\|_{B^\alpha_{q ,c(\mathbb{N})}(\mathbb{T}^3)} .
\ea\ee
\end{lemma}
\begin{proof}
Firstly, by direct calculation, we have
\be\label{2.15}\ba
\nabla\left(\frac{(\rho v)^\varepsilon}{\rho^\varepsilon}\right)
&=
\frac{\nabla(\rho v)^\varepsilon - v\,\nabla\rho^\varepsilon}{\rho^\varepsilon}
-
\frac{(\rho v)^\varepsilon-\rho^\varepsilon v}{(\rho^\varepsilon)^2}\,
\nabla\rho^\varepsilon
=: I_1 + I_2 .
\ea\ee
Making use of  the Minkowski inequality and  Corollary \ref{coro}, we obtain
\be\label{2.16}\ba
\|I_1\|_{L^q_x}
&\le
c_1^{-1}
\int_{B_\varepsilon(0)}
|\nabla\varrho_\varepsilon(y)|\,
\|\rho(x-y)\|_{L^\infty_x}\|v(x-y)-v(x)\|_{L^q_x}\,dy \\
&\leq C\int_{B_\varepsilon(0)}\f{1}{\varepsilon^4}|\nabla \varrho(\f{y}{\varepsilon})|\|\rho\|_{L^\infty_x}\frac{\|v(x-y)-v(x)\|_{L^q_x}}{|y|^\alpha}|y|^\alpha\, dy\\
&\le
C\varepsilon^{\alpha-1}\|\rho\|_{L^\infty_x}\,\|v\|_{B^\alpha_{q,c(\mathbb{N})}}
\int_{B_\varepsilon(0)}\f{1}{\varepsilon^3}
|\nabla\varrho(\f{y}{\varepsilon})|\,dy\\
&\leq C\varepsilon^{\alpha-1}\|\rho\|_{L^\infty_x}\,\|v\|_{B^\alpha_{q,c(\mathbb{N})}}.
\ea\ee
Regarding term $I_2$, it follows from the  H\"older inequality that
\be\label{2.17}\ba
\|I_2\|_{L^q_x}
&\le
c_1^{-2}
\left\|
(\rho v)^\varepsilon-\rho^\varepsilon v\right\|_{L^q_x}\|\nabla \rho^\varepsilon\|_{L^\infty_x}.
\ea\ee
Note that
\begin{align*}
(\rho v)^\varepsilon(x)-\rho^\varepsilon(x)v(x)
=
\int_{\mathbb{T}^3}\varrho_\varepsilon(y)\,
\rho(x-y)\bigl(v(x-y)-v(x)\bigr)\,dy.
\end{align*}
Following the same path  of \eqref{2.16}, we obtain
\be\label{2.18}
\left\|
(\rho v)^\varepsilon-\rho^\varepsilon v
\right\|_{L^q_x}
\le
C\,\varepsilon^\alpha\,
\|\rho\|_{L^\infty}\,
\|v\|_{B^\alpha_{q,c(\mathbb{N})}}.
\ee
On the other hand, using $\nabla\rho^\varepsilon=\rho *\nabla\varrho_\varepsilon$
and Young's inequality, we have
\begin{align*}
\|\nabla\rho^\varepsilon\|_{L^\infty}
\le
\|\rho\|_{L^\infty}\,
\|\nabla\varrho_\varepsilon\|_{L^1}
\le
C\,\varepsilon^{-1}\,
\|\rho\|_{L^\infty},
\end{align*}
from which it follows that
\be\label{2.19}\ba
\|I_2\|_{L^q_x}
\le
C\,\varepsilon^{\alpha-1}\,
\|\rho\|_{L^\infty}\,
\|v\|_{B^\alpha_{q,c(\mathbb{N})}}.
\ea\ee
Combining  \eqref{2.17} and \eqref{2.19}, we conclude that
\begin{align*}
\left\|
\nabla\left(\frac{(\rho v)^\varepsilon}{\rho^\varepsilon}\right)
\right\|_{L^q_x}
\le
C\,\varepsilon^{\alpha-1}\,
\|\rho\|_{L^\infty}\,
\|v\|_{B^\alpha_{q,c(\mathbb{N})}}.
\end{align*}
Thus we have proved the desired estimate \eqref{2.13}.

Next, we are in a position to prove the estimate \eqref{2.14}.  Due to the Leibniz's rule, we have
\be\label{2.20}\ba
&\nabla^2\!\left(\frac{(\rho v)^\varepsilon}{\rho^\varepsilon}\right)\\
=&
\frac{\rho^\varepsilon\nabla^2(\rho v)^\varepsilon-(\rho v)^\varepsilon\nabla^2\rho^\varepsilon}{(\rho^\varepsilon)^2}
-\frac{\nabla \frac{(\rho v)^\varepsilon}{\rho^\varepsilon}\otimes \nabla \rho^\varepsilon+\nabla \rho^\varepsilon\otimes\nabla \frac{(\rho v)^\varepsilon}{\rho^\varepsilon}}{\rho^\varepsilon} \\
=&
\frac{\nabla^2(\rho v)^\varepsilon-v\,\nabla^2\rho^\varepsilon}{\rho^\varepsilon}
+\frac{\rho^\varepsilon v-(\rho v)^\varepsilon}{(\rho^\varepsilon)^2}\nabla^2\rho^\varepsilon
-\frac{\nabla \frac{(\rho v)^\varepsilon}{\rho^\varepsilon}\otimes \nabla \rho^\varepsilon+\nabla \rho^\varepsilon\otimes\nabla \frac{(\rho v)^\varepsilon}{\rho^\varepsilon}}{\rho^\varepsilon} \\
=&: J_1+J_2+J_3.
\ea\ee
For term $J_1$, by the same token as \eqref{2.16}, we infer that
\be\label{2.21}\ba
\|J_1\|_{L^q_x}
&\le
c_1^{-2}
\int_{B_\varepsilon(0)}
|\nabla^2\varrho_\varepsilon(y)|\,
\|\rho(x-y)\|_{L^\infty_x}\|v(x-y)-v(x)\|_{L^q_x}\,dy \\
&\leq C\int_{B_\varepsilon(0)}\f{1}{\varepsilon^5}|\nabla^2 \varrho(\f{y}{\varepsilon})|\|\rho\|_{L^\infty_x}\frac{\|v(x-y)-v(x)\|_{L^q_x}}{|y|^\alpha}|y|^\alpha\, dy\\
&\le
C\varepsilon^{\alpha-2}\|\rho\|_{L^\infty}\,\|v\|_{B^\alpha_{q,c(\mathbb{N})}}
\int_{B_\varepsilon(0)}\f{1}{\varepsilon^3}
|\nabla^2\varrho(\f{y}{\varepsilon})|\,dy\\
&\leq C\varepsilon^{\alpha-2}\|\rho\|_{L^\infty}\,\|v\|_{B^\alpha_{q,c(\mathbb{N})}}.
\ea\ee
By using the H\"older inequality, \eqref{2.18} and standard mollifier property, we know that
\be\label{2.22}\ba
\|J_2\|_{L^q_x}&\leq c_1^{-2}\left\|
(\rho v)^\varepsilon-\rho^\varepsilon v
\right\|_{L^q_x}\|\nabla ^2\rho^\varepsilon \|_{L^\infty_x}\\
&\leq C\varepsilon^\alpha\,
\|\rho\|_{L^\infty_x}\,
\|v\|_{B^\alpha_{q,c(\mathbb{N})}}\varepsilon^{-2}\|\rho\|_{L^\infty_x}\\
&\leq C\varepsilon^{\alpha-2}\|\rho\|_{L^\infty_x}^2\|v\|_{B^\alpha_{q,c(\mathbb{N})}}.
\ea\ee
For the last term $J_3$, in light of the H\"older inequality and \eqref{2.13}, we have
\be\label{2.23}\ba
\|J_3\|_{L^q_x}&\leq c_1^{-1}\|\nabla \left(\f{(\rho v)^\varepsilon}{\rho^\varepsilon}\right)\|_{L^q_x}\|\nabla \rho^\varepsilon\|_{L^\infty_x}\\
&\leq C\,\varepsilon^{\alpha-1}\,
\|v\|_{B^\alpha_{q,c(\mathbb{N})}} \varepsilon^{-1}\|\rho\|_{L^\infty_x}\\
&\leq C\,\varepsilon^{\alpha-2}\,
\|v\|_{B^\alpha_{q,c(\mathbb{N})}}.
\ea\ee
Then combining \eqref{2.21}, \eqref{2.22} and \eqref{2.23} yields
\eqref{2.14}, and the proof of this lemma  is completed.
\end{proof}
To deal with conditions (iii) and (iv) of Theorem \ref{the 1.1}, we need the following Sobolev version of Lemma \ref{lem2.4}, whose part (i) was proved in \cite[Lemma 2.3]{[WY]}.
\begin{lemma}\label{lem-sob}
Let $1\leq p,q<\infty$, $0<\beta<1$ and $0<c_1\leq \rho\leq c_2<\infty$.
\begin{enumerate}[(i)]
\item If $v\in W^{1,q}(\mathbb{T}^3)$, then there exists a constant $C>0$ independent of $\varepsilon$ such that
\be\label{2.24}
\Big\|\nabla\Big(\frac{(\rho v)^\varepsilon}{\rho^\varepsilon}\Big)\Big\|_{L^q(\mathbb{T}^3)}\leq C\|\nabla v\|_{L^q(\mathbb{T}^3)}.
\ee
\item If $\rho\in L^\infty(0,T;B^{\beta}_{\infty,\infty}(\mathbb{T}^3))$ and $v\in L^p(0,T;W^{1,q}(\mathbb{T}^3))$, then
\be\label{2.25}
\lim_{\varepsilon\to0}\varepsilon\Big\|\nabla^2\Big(\frac{(\rho v)^\varepsilon}{\rho^\varepsilon}\Big)\Big\|_{L^p(0,T;L^q(\mathbb{T}^3))}=0.
\ee
\end{enumerate}
\end{lemma}
\begin{proof}
	Estimate \eqref{2.24} follows from \cite[Lemma 2.3]{[WY]}.
	It remains to prove \eqref{2.25}.
	Throughout the proof, $C>0$ is independent of $\varepsilon$ and may
	depend on $c_1$, $c_2$, the fixed exponents, the mollifier and
	$\|\rho\|_{L^\infty(0,T;B^\beta_{\infty,\infty}(\mathbb{T}^3))}$.
	
   Before going any further, we conclude from a direct calculation and  Lemma~\ref{lem2.1} that 
       \be\label{2.26}\ba
    \|(\rho v)^\varepsilon-\rho^\varepsilon v\|_{L^q}
    &\leq c_2\int_{B_\varepsilon(0)}\varrho_\varepsilon(y)
    \|v(\cdot-y)-v(\cdot)\|_{L^q}\,dy\\
    &\leq c_2\|\nabla v\|_{L^q}
    \int_{B_\varepsilon(0)}\varrho_\varepsilon(y)|y|\,dy
    \leq c_2\,\varepsilon\|\nabla v\|_{L^q},\\
    \|\nabla^k\rho^\varepsilon\|_{L^\infty}
    &\leq C\varepsilon^{\beta-k}
    \|\rho\|_{B^\beta_{\infty,\infty}(\mathbb{T}^3)}
    \leq C\varepsilon^{\beta-k},\quad k=1,2.
  \ea\ee
	
Recall \eqref{2.20}, we estimate $J_1$, $J_2$ and $J_3$ as follows.
	For $J_1$, using
	$\int_{\mathbb{R}^3}\nabla^2\varrho_\varepsilon(y)\,dy=0$,
	we have
	\be\label{2.27}\ba
	\nabla^2(\rho v)^\varepsilon-v\nabla^2\rho^\varepsilon
	&=\int_{B_\varepsilon(0)}
	\nabla^2\varrho_\varepsilon(y)\rho(x-y)
	\big(v(x-y)-v(x)\big)\,dy\\
	&=\int_{B_\varepsilon(0)}\nabla^2\varrho_\varepsilon(y)
	\big[\rho(x)+\big(\rho(x-y)-\rho(x)\big)\big]
	\big(v(x-y)-v(x)\big)\,dy\\
	&=\rho\,\nabla^2v^\varepsilon
	+\int_{B_\varepsilon(0)}
	\nabla^2\varrho_\varepsilon(y)
	\big(\rho(x-y)-\rho(x)\big)
	\big(v(x-y)-v(x)\big)\,dy.
	\ea\ee
	Applying Minkowski's inequality, using \eqref{2.2} and the Sobolev translation estimate, we obtain
	\be\label{2.28}\ba
	\|J_1\|_{L^q}
	&\leq c_1^{-1}\Big(
	c_2\|\nabla^2v^\varepsilon\|_{L^q}
	+\int_{B_\varepsilon(0)}
	|\nabla^2\varrho_\varepsilon(y)|\,
	\|\rho(\cdot-y)-\rho(\cdot)\|_{L^\infty}
	\|v(\cdot-y)-v(\cdot)\|_{L^q}\,dy
	\Big)\\
	&\leq \frac{c_2}{c_1}\|\nabla^2v^\varepsilon\|_{L^q}
	+C\varepsilon^{\beta-1}\|\nabla v\|_{L^q}.
	\ea\ee
	For $J_2$ and $J_3$, applying H\"older's inequality
	and using \eqref{2.26} and \eqref{2.24}, we obtain
	\be\label{2.29}\ba
	\|J_2\|_{L^q}
	&\leq c_1^{-2}
	\|(\rho v)^\varepsilon-\rho^\varepsilon v\|_{L^q}
	\|\nabla^2\rho^\varepsilon\|_{L^\infty}
	\leq C\varepsilon^{\beta-1}\|\nabla v\|_{L^q},\\
	\|J_3\|_{L^q}
	&\leq 2c_1^{-1}
	\Big\|\nabla\Big(\frac{(\rho v)^\varepsilon}
	{\rho^\varepsilon}\Big)\Big\|_{L^q}
	\|\nabla\rho^\varepsilon\|_{L^\infty}
	\leq C\varepsilon^{\beta-1}\|\nabla v\|_{L^q}.
	\ea\ee
	Combining \eqref{2.20}, \eqref{2.28} and \eqref{2.29},
	we obtain, for a.e. $t\in(0,T)$,
	\be\label{2.30}
	\varepsilon
	\Big\|\nabla^2\Big(\frac{(\rho v)^\varepsilon}
	{\rho^\varepsilon}\Big)(t)\Big\|_{L^q}
	\leq \frac{c_2}{c_1}\,
	\varepsilon\|\nabla^2v^\varepsilon(t)\|_{L^q}
	+C\varepsilon^\beta\|\nabla v(t)\|_{L^q}.
	\ee
	
    Taking the $L^p(0,T)$ norm in \eqref{2.30} and applying
    \cite[Lemma 2.1]{[WY]} to $f=\nabla v$ with $p,q<\infty$,
    we obtain
	\begin{equation*}\ba
		\varepsilon
		\Big\|\nabla^2\Big(\frac{(\rho v)^\varepsilon}
		{\rho^\varepsilon}\Big)\Big\|_{L^p(0,T;L^q)}
		&\leq \frac{c_2}{c_1}\,
		\varepsilon\|\nabla^2v^\varepsilon\|_{L^p(0,T;L^q)}
		+C\varepsilon^\beta
		\|\nabla v\|_{L^p(0,T;L^q)}\to0\qquad(\varepsilon\to0),
		\ea\end{equation*}
	which proves \eqref{2.25}.
\end{proof}

To prove the endpoint case (iv) of Theorem \ref{the 1.1}, we need the following lemma, whose original version can be found in \cite[Lemma 7.7]{[GT]}.
\begin{lemma} \label{lem-zero-set}
Let $1\leq p\leq\infty$ and $f\in W^{1,p}(\mathbb{T}^3)$. Then
\begin{equation*}
	\nabla f=0\quad\text{a.e. on }\{x\in\mathbb{T}^3:f(x)=0\}.
\end{equation*}
\end{lemma}
\begin{proof}
	Since the proof of this lemma is analogous to that of Lemma 7.7 in \cite{[GT]}, we omit the details.
\end{proof}

Eventually, for the convenience of readers, we end this section by presenting the definition of weak solutions to the compressible Euler equations \eqref{CE}.

 \begin{definition}\label{ceulerdefi}
 	A pair ($\rho,v$) is called a weak solution to  the   compressible Euler equations  \eqref{CE}  if ($\rho,v$) satisfies
 	 	\begin{enumerate}[(i)]
 	 		\item $\sqrt{\rho}v\in L^\infty(0,T;L^{2}(\Omega))$  and $\rho\int_1^\rho \frac{\pi(s)}{s^2 } ds\in  L^\infty(0,T;L^{1}(\Omega))$;
 		\item For any test function $\Phi\in C_0^\infty((0,T)\times \Omega)$, there holds 		\begin{equation}\label{2.31}\int^T_0\int_{\Omega}\rho(t,x)\partial_t \Phi(t,x)+ \rho(t,x)v(t,x)\nabla \Phi(t,x)dxdt=0;\end{equation}
 		\item
 		For any  test vector field $\Psi\in C_{0}^{\infty}((0,T)\times\Omega)$, there holds
 		\begin{equation}\label{2.32}
 			\begin{aligned} 				\int_{0}^{T}\int_{\Omega}\rho(t,x)v(t,x)&\partial_{t}\Psi(t,x)+\B(\rho(t,x)v(t,x)\otimes v(t,x)\B) \nabla\Psi(t,x)\\
 				+ &
 				\pi(\rho )(t,x)\text{div} \Psi(t,x)dxdt=0.
 		\end{aligned}	\end{equation}
 		 	\end{enumerate}
 \end{definition}

\section{Helicity  conservation for the compressible Euler equations}
In this section, we are devoted to the proof of our main result regarding the helicity conservation of weak solutions to the compressible Euler equations. We will divide the proof into two steps. Precisely, we first show the local helicity conservation in the sense of distributions of time, next we extend the local  helicity identity to be the global one.
\begin{proof}[Proof of Theorem~\ref{the 1.1}.]
	\underline{Step 1:} First, we show that the helicity of the weak solutions to the compressible Euler equations is locally conserved in the sense of distribution of $(0,T)$. To this end, let $\phi(t)$ be a smooth function compactly supported in $(0,T)$, multiplying $\eqref{CE}_2$ by  $\varphi^\varepsilon
	:=
	\left(\phi(t)
	\frac{1}{\rho^\varepsilon}\,
	\nabla\times
	\frac{(\rho v)^\varepsilon}{\rho^\varepsilon}
	\right)^\varepsilon $ and integrating the resultant equation over $(0,T)\times \mathbb{T}^3$, we have
\begin{align}\label{3.1}
\int_0^T\int_{\mathbb{T}^3}
\phi(t)\frac{1}{\rho^\varepsilon}\,
\nabla\times
\frac{(\rho v)^\varepsilon}{\rho^\varepsilon}\Big[
(\rho v)_t^\varepsilon+\s(\rho v\otimes v)^\varepsilon+\nabla\pi(\rho)^\varepsilon
\Big]\,dxdt
=0.
\end{align}
For the sake of convenience, we denote the terms on the left-hand side of \eqref{3.1} by $L_1-L_3$, respectively. Regarding the term $L_1$, a straightforward calculation yields that
\be\label{3.2}\begin{aligned}
L_1
&=\int_0^T\int_{\mathbb T^3}\phi(t)\frac{(\rho v)_t^\varepsilon}{\rho^\varepsilon}\cdot
\nabla\times\Big(\frac{(\rho v)^\varepsilon}{\rho^\varepsilon}\Big)\,dxdt \\
&=\int_0^T\int_{\mathbb T^3}\phi(t)\B[
\partial_t\!\left(\frac{(\rho v)^\varepsilon}{\rho^\varepsilon}\right)+\frac{(\rho v)^\varepsilon\,\rho_t^\varepsilon}{(\rho^\varepsilon)^2}\B]\cdot
\nabla\times\left(\frac{(\rho v)^\varepsilon}{\rho^\varepsilon}\right)\,dxdt \\
&=-\frac12\int_0^T \partial_t \phi(t)\int_{\mathbb{T}^3}
\frac{(\rho v)^\varepsilon}{\rho^\varepsilon}\cdot
\nabla\times\left(\frac{(\rho v)^\varepsilon}{\rho^\varepsilon}\right)\,dxdt \\
&\,\,\,\,\,\,\,\,\,-\int_0^T\phi(t)\int_{\mathbb{T}^3} \frac{(\rho v)^\varepsilon}{(\rho^\varepsilon)^2}\s (\rho v)^\varepsilon \cdot \nabla \times \left(\frac{(\rho v)^\varepsilon}{\rho^\varepsilon}\right)\, dxdt\\
&=L_{11}+L_{12},
\end{aligned}\ee
where we have used the integration by parts, the equation $\eqref{CE}_{1}$ and the identity \eqref{eqrateofhelicity}.
Then, in light of the integration by parts again, we reformulate  the term $L_{11}$  as
\be\label{3.3}\ba
L_{11}&=-\f{1}{2}\int_0^T \partial_t \phi(t) \int_{\mathbb T^3} \left(v^\varepsilon
+
\frac{(\rho v)^\varepsilon-\rho^\varepsilon v^\varepsilon}{\rho^\varepsilon}\right)\cdot \nabla \times \left(v^\varepsilon
+
\frac{(\rho v)^\varepsilon-\rho^\varepsilon v^\varepsilon}{\rho^\varepsilon}\right)\, dxdt\\
&=-\f{1}{2}\int_0^T \partial_t \phi(t) \B[\int_{\mathbb{T}^3} v^\varepsilon\cdot(\nabla\times v^\varepsilon)\,dx+ \int_{\mathbb{T}^3}
v^\varepsilon\cdot
\nabla\times\left(\frac{(\rho v)^\varepsilon-\rho^\varepsilon v^\varepsilon}{\rho^\varepsilon}\right)\,dx\\
&\quad+
\int_{\mathbb{T}^3}
\frac{(\rho v)^\varepsilon-\rho^\varepsilon v^\varepsilon}{\rho^\varepsilon}\cdot
\nabla\times\left(\frac{(\rho v)^\varepsilon}{\rho^\varepsilon}\right)\,dx
\B] dt\\
&=-\f{1}{2}\int_0^T\int_{\mathbb T^3} v^\varepsilon \cdot \omega^\varepsilon \partial_t \phi(t) \, dxdt\\
&\quad-\f{1}{2} \int_0^T\int_{\mathbb T^3} \frac{(\rho v)^\varepsilon-\rho^\varepsilon v^\varepsilon}{\rho^\varepsilon}\cdot\B[\omega^\varepsilon +\nabla \times \left(\frac{(\rho v)^\varepsilon}{\rho^\varepsilon}\right) \B]\partial_t \phi(t)\, dxdt.
\ea\ee
By the same token, we rewrite term $L_{12} $ as
\be\label{3.4}\ba
L_{12}&=-\int_0^T\phi(t)\int_{\mathbb{T}^3}\left( \frac{(\rho v)^\varepsilon-\rho^\varepsilon v^\varepsilon}{(\rho^\varepsilon)^2}+\f{v^\varepsilon}{\rho ^\varepsilon}\right)\s (\rho v)^\varepsilon \cdot \nabla \times \left(\frac{(\rho v)^\varepsilon}{\rho^\varepsilon}\right)\, dxdt\\
&=-\int_0^T\phi(t)\int_{\mathbb{T}^3} \frac{(\rho v)^\varepsilon-\rho^\varepsilon v^\varepsilon}{(\rho^\varepsilon)^2}\s (\rho v)^\varepsilon \cdot \nabla \times \left(\frac{(\rho v)^\varepsilon}{\rho^\varepsilon}\right)\, dxdt\\
&\quad-\int_0^T\phi(t)\int_{\mathbb{T}^3} \f{v^\varepsilon}{\rho ^\varepsilon}\s (\rho v)^\varepsilon \cdot \nabla \times \left(\frac{(\rho v)^\varepsilon}{\rho^\varepsilon}\right)\, dxdt.
\ea\ee
Substituting \eqref{3.3} and \eqref{3.4} into \eqref{3.2}, we obtain
\be\label{3.5}\ba
L_1&=-\f{1}{2}\int_0^T\int_{\mathbb T^3} v^\varepsilon \cdot \omega^\varepsilon \partial_t \phi(t) \, dxdt\\
&\quad-\f{1}{2} \int_0^T\int_{\mathbb T^3} \frac{(\rho v)^\varepsilon-\rho^\varepsilon v^\varepsilon}{\rho^\varepsilon}\cdot\B[\omega^\varepsilon +\nabla \times \left(\frac{(\rho v)^\varepsilon}{\rho^\varepsilon}\right) \B]\partial_t \phi(t)\, dxdt\\
&\quad-\int_0^T\int_{\mathbb{T}^3} \phi(t)\frac{(\rho v)^\varepsilon-\rho^\varepsilon v^\varepsilon}{(\rho^\varepsilon)^2}\s (\rho v)^\varepsilon \cdot \nabla \times \left(\frac{(\rho v)^\varepsilon}{\rho^\varepsilon}\right)\, dxdt\\
&\quad-\int_0^T\int_{\mathbb{T}^3}\phi(t) \f{v^\varepsilon}{\rho ^\varepsilon}\s (\rho v)^\varepsilon \cdot \nabla \times \left(\frac{(\rho v)^\varepsilon}{\rho^\varepsilon}\right)\, dxdt.
\ea\ee
Regarding the term $L_2$, it follows from the integration by parts that
\be\label{3.6}\ba
L_2&=\int_0^T \int_{\mathbb{T}^3}\phi(t)\frac{1}{\rho^\varepsilon}\,
\nabla\times\Big(\frac{(\rho v)^\varepsilon}{\rho^\varepsilon}\Big)\cdot
\s\Big[(\rho v\otimes v)^\varepsilon-(\rho v)^\varepsilon\otimes v^\varepsilon\Big]\,dxdt \\
&\quad
+\int_0^T\int_{\mathbb{T}^3}\phi(t)\frac{1}{\rho^\varepsilon}\,
\nabla\times\Big(\frac{(\rho v)^\varepsilon}{\rho^\varepsilon}\Big)\cdot
\s\big[(\rho v)^\varepsilon\otimes v^\varepsilon\big]\,dxdt\\
&=-
\int_0^T\int_{\mathbb{T}^3}\phi(t)
\Big[(\rho v\otimes v)^\varepsilon-(\rho v)^\varepsilon\otimes v^\varepsilon\Big]:
\nabla\!\left(
\frac{1}{\rho^\varepsilon}\,
\nabla\times\Big(\frac{(\rho v)^\varepsilon}{\rho^\varepsilon}\Big)
\right)\,dxdt\\
&\quad
+\int_0^T\int_{\mathbb{T}^3}\phi(t)
\frac{v^\varepsilon}{\rho^\varepsilon}\,
\s(\rho v)^\varepsilon\cdot
\nabla\times\Big(\frac{(\rho v)^\varepsilon}{\rho^\varepsilon}\Big)\,dxdt\\
&\quad+\int_0^T\int_{\mathbb{T}^3}\phi(t)
\frac{1}{\rho^\varepsilon}\,
\B((\rho v)^\varepsilon\cdot\nabla v^\varepsilon\B) \cdot
\nabla\times\Big(\frac{(\rho v)^\varepsilon}{\rho^\varepsilon}\Big)\,dxdt.
\ea\ee
Next, we deal with the pressure term $L_3$, by virtue of the vector identity $\nabla\cdot(\nabla\times F)=0$ and integration  by parts, we know that
\be\label{3.7}\ba
L_3
&=\int_0^T\int_{\mathbb{T}^3}\phi(t)\frac{1}{\rho^\varepsilon}\,
\nabla\times\Big(\frac{(\rho v)^\varepsilon}{\rho^\varepsilon}\Big)\cdot
\big(\nabla(\pi(\rho))^\varepsilon-\nabla\pi(\rho^\varepsilon)\big)\,dxdt\\
&\quad+\int_0^T\int_{\mathbb{T}^3}\phi(t)\frac{1}{\rho^\varepsilon}\,
\nabla\times\Big(\frac{(\rho v)^\varepsilon}{\rho^\varepsilon}\Big)\cdot \nabla\pi(\rho^\varepsilon)\,dx dt\\
&=\int_0^T\int_{\mathbb{T}^3}\phi(t)\frac{1}{\rho^\varepsilon}\,
\nabla\times\Big(\frac{(\rho v)^\varepsilon}{\rho^\varepsilon}\Big)\cdot
\big(\nabla(\pi(\rho))^\varepsilon-\nabla\pi(\rho^\varepsilon)\big)\,dxdt\\
&\quad+\int_0^T\int_{\mathbb{T}^3}\phi(t)
\nabla\Big(\int_{1}^{\rho^\varepsilon}\frac{\pi'(s)}{s}\,ds\Big)\cdot
\nabla\times\Big(\frac{(\rho v)^\varepsilon}{\rho^\varepsilon}\Big)\,dxdt\\
&=\int_0^T\int_{\mathbb{T}^3}\phi(t)\frac{1}{\rho^\varepsilon}\,
\nabla\times\Big(\frac{(\rho v)^\varepsilon}{\rho^\varepsilon}\Big)\cdot
\big(\nabla(\pi(\rho))^\varepsilon-\nabla\pi(\rho^\varepsilon)\big)\,dxdt.
\ea\ee
Now, plugging \eqref{3.5}-\eqref{3.7} into \eqref{3.1}, we have
\be\label{3.8}\ba
&\f{1}{2}\int_0^T\int_{\mathbb T^3} v^\varepsilon \cdot \omega^\varepsilon \partial_t \phi(t) \, dxdt\\
=&-\int_0^T\int_{\mathbb{T}^3}\phi(t)
\Big[(\rho v\otimes v)^\varepsilon-(\rho v)^\varepsilon\otimes v^\varepsilon\Big]:
\nabla\!\left(
\frac{1}{\rho^\varepsilon}\,
\nabla\times\Big(\frac{(\rho v)^\varepsilon}{\rho^\varepsilon}\Big)
\right)\,dxdt\\
&+\int_0^T\int_{\mathbb{T}^3}\phi(t)
\frac{1}{\rho^\varepsilon}\,
\B((\rho v)^\varepsilon\cdot\nabla v^\varepsilon\B) \cdot
\nabla\times\Big(\frac{(\rho v)^\varepsilon}{\rho^\varepsilon}\Big)\,dxdt\\
&-\f{1}{2} \int_0^T\int_{\mathbb T^3} \frac{(\rho v)^\varepsilon-\rho^\varepsilon v^\varepsilon}{\rho^\varepsilon}\cdot\B[\omega^\varepsilon +\nabla \times \left(\frac{(\rho v)^\varepsilon}{\rho^\varepsilon}\right) \B]\partial_t \phi(t)\, dxdt\\
&-\int_0^T\int_{\mathbb{T}^3} \phi(t)\frac{(\rho v)^\varepsilon-\rho^\varepsilon v^\varepsilon}{(\rho^\varepsilon)^2}\s (\rho v)^\varepsilon \cdot \nabla \times \left(\frac{(\rho v)^\varepsilon}{\rho^\varepsilon}\right)\, dxdt\\
&+\int_0^T\int_{\mathbb{T}^3}\phi(t)\frac{1}{\rho^\varepsilon}\,
\nabla\times\Big(\frac{(\rho v)^\varepsilon}{\rho^\varepsilon}\Big)\cdot
\big(\nabla(\pi(\rho))^\varepsilon-\nabla\pi(\rho^\varepsilon)\big)\,dxdt\\
=&
R_1+R_2+R_3+R_4+R_5.
\ea\ee
To obtain the desired local helicity equality, it suffices to show the terms $R_1-R_5$ tend to zero as $\varepsilon \rightarrow 0$. We  will complete the proof by considering three subcases.

{\bf General case: $\f23<\alpha <1$.}

Before going further, we presuppose that 
\begin{equation*}
	\begin{aligned}
	&	v,\rho v\in L^{\theta p_1}(0,T;B^{\alpha}_{\theta q_1, c(\mathbb{N})}(\mathbb{T}^3))\cap L^\infty(0,T;L^2(\mathbb{T}^3)),\\
	& 0<c_1<\rho<c_2<\infty,~ \rho \in L^{\gamma p_3}(0,T;B^{\beta}_{\gamma q_3, \infty}(\mathbb{T}^3)),	
	\end{aligned}
\end{equation*} 
which will be determined later. Firstly, by means of  H\"older's inequality, we arrive at
\be\label{3.9}\ba |R_1|&=\B|\int_0^T\int_{\mathbb{T}^3}\phi(t)
\Big[(\rho v\otimes v)^\varepsilon-(\rho v)^\varepsilon\otimes v^\varepsilon\Big]:
\B[\!\left(\nabla\times(\frac{(\rho v)^\varepsilon}{\rho^\varepsilon})\right)\otimes
\nabla\!(\frac{1}{\rho^\varepsilon})\\
&\quad
+
\frac{1}{\rho^\varepsilon}\,
\nabla\!\left(\nabla\times(\frac{(\rho v)^\varepsilon}{\rho^\varepsilon})\right)\B]\,dxdt\B|\\
&\le
C\,
\|(\rho v\otimes v)^\varepsilon-(\rho v)^\varepsilon\otimes v^\varepsilon\|_{L^{\frac{p}{p-1}}_t(L^{\frac{q}{q-1}}_x)}  \\
&\quad
\Big(
\|(\nabla\times(\frac{(\rho v)^\varepsilon}{\rho^\varepsilon}))\otimes\nabla(\tfrac1{\rho^\varepsilon})\|_{L^p_t(L^q_x)}
+
\|\tfrac1{\rho^\varepsilon}\nabla(\nabla\times\tfrac{(\rho v)^\varepsilon}{\rho^\varepsilon})\|_{L^p_t(L^q_x)}
\Big).
\ea\ee

On the one hand, applying Lemma \ref{lem2.3} by choosing $f= \rho v, g= v, r=\f{p}{p-1}, s=\f{q}{q-1}, \theta p_1= \gamma p_3=p,\theta q_1=\gamma q_3=q,$ and $ (1-\theta)p_2=(1-\gamma) p_4=\infty, (1-\theta )q_2=(1-\gamma )q_4=2$, we have 
\be\label{3.10}
\ba
\|(\rho v\otimes v)^\varepsilon-(\rho v)^\varepsilon\otimes v^\varepsilon\|_{L^{\frac{p}{p-1}}_t(L^{\frac{q}{q-1}}_x)} \leq o(\varepsilon ^{2\theta \alpha}).
\ea
\ee
On the other hand, taking advantage of Lemma \ref{lem2.4}, we obtain
\be\label{3.11}
\ba
&\|(\nabla\times(\frac{(\rho v)^\varepsilon}{\rho^\varepsilon}))\otimes\nabla(\tfrac1{\rho^\varepsilon})\|_{L^p_t(L^q_x)}\\\leq &C\|\nabla \times (\frac{(\rho v)^\varepsilon}{\rho^\varepsilon})\|_{L^p_t(L^q_x)}\|\nabla \f{1}{\rho^\varepsilon}\|_{L^\infty_t(L^\infty_x)}\\
\leq& o(\varepsilon ^{\alpha -1})O(\varepsilon^{-1})\leq o(\varepsilon^{\alpha-2}),
\ea\ee
and
\be\label{3.12}\ba
\|\tfrac1{\rho^\varepsilon}\nabla(\nabla\times\tfrac{(\rho v)^\varepsilon}{\rho^\varepsilon})\|_{L^p_t(L^q_x)}\leq o(\varepsilon^{\alpha -2}).
\ea\ee
Then, combining \eqref{3.10}-\eqref{3.12}, we conclude that 
\be\label{3.13}
\ba
|R_1|\leq o(\varepsilon^{2\theta \alpha +\alpha -2}).
\ea\ee
Hence, letting $2\theta \alpha +\alpha -2=0$ with $\f{2}{3}<\alpha <1$, which implies $\rho v, v\in L^{\f{2}{\alpha}}(0,T;B^{\alpha}_{\f{4}{2-\alpha},c(\mathbb{N})}(\mathbb{T}^3))\cap L^\infty(0,T;L^2(\mathbb{T}^3))$ ensures term $R_1$ converges to zero as $\varepsilon \to 0$. 

Next, by direct computation, the second term $R_2$ can be rewritten as
\be\label{3.14}\ba
|R_2|
&\le 
\B|\int_0^T \int_{\mathbb{T}^3}\phi(t)
\Big(\nabla\times\frac{(\rho v)^\varepsilon}{\rho^\varepsilon}\Big)\cdot
\big(\nabla v^\varepsilon\cdot v^\varepsilon\big)
\,dxdt\B|
\\
&\quad+\B|
\int_0^T \int_{\mathbb{T}^3} \phi(t)
\Big(\nabla\times\frac{(\rho v)^\varepsilon}{\rho^\varepsilon}\Big)\cdot
\Big(\nabla v^\varepsilon\cdot \frac{(\rho v)^\varepsilon-\rho^\varepsilon v^\varepsilon}{\rho^\varepsilon}\Big)
\,dxdt\B|\\
&=R_{21}+R_{22}.
\ea\ee
Regarding term $R_{21}$, by using the vector identity
$
(v^\varepsilon\cdot\nabla)v^\varepsilon=\frac12\nabla|v^\varepsilon|^2+(\nabla\times v^\varepsilon)\times v^\varepsilon,
$
we arrive at 
\be\label{3.15}\ba
|R_{21}|
&=\B|
\int_0^T \int_{\mathbb{T}^3}\frac12\phi(t)
\nabla|v^\varepsilon|^2\cdot
\nabla\times\Big(\frac{(\rho v)^\varepsilon}{\rho^\varepsilon}\Big)
\,dxdt
\\
&\quad+
\int_0^T \int_{\mathbb{T}^3}\phi(t)
(\nabla\times v^\varepsilon)\times v^\varepsilon
\cdot
\nabla\times\Big(\frac{(\rho v)^\varepsilon}{\rho^\varepsilon}\Big)
\,dxdt
\B|\\
&= \B|\int_0^T \int_{\mathbb{T}^3}\phi(t)
(\omega^\varepsilon\times v^\varepsilon)
\cdot\left[
\nabla\times\Big(\frac{(\rho v)^\varepsilon-\rho^\varepsilon v^\varepsilon}{\rho^\varepsilon}\Big)+\omega^\varepsilon\right]
\,dxdt\B|\\
&= \B|\int_0^T \int_{\mathbb{T}^3}\phi(t)
\nabla\times(\omega^\varepsilon\times v^\varepsilon)
\cdot
\Big(\frac{(\rho v)^\varepsilon-\rho^\varepsilon v^\varepsilon}{\rho^\varepsilon}\Big)
\,dxdt\B|,\\
\ea\ee
where the integration by parts,  the vector identities $\s (\nabla \times {F})=0$ and $\omega^\varepsilon \times v^\varepsilon \cdot \omega^\varepsilon=\omega^\varepsilon \times \omega^\varepsilon  \cdot v^\varepsilon=0$ have been used.

To deal with the term on the right-hand side of \eqref{3.15}, combining the   vector identity
$
\nabla\times(w^\varepsilon\times v^\varepsilon)=w^\varepsilon\s v^\varepsilon - v^\varepsilon\s w^\varepsilon + v^\varepsilon\cdot\nabla w^\varepsilon - w^\varepsilon\cdot\nabla v^\varepsilon,
$
and  $\s w^\varepsilon=0$, we can further infer that 
\be\label{3.16}\ba
|R_{21}|&\leq \B|\int_0^T\int_{\mathbb{T}^3} \phi(t)\left(w^\varepsilon\s v^\varepsilon  + v^\varepsilon\cdot\nabla w^\varepsilon - w^\varepsilon\cdot\nabla v^\varepsilon\right)\cdot \Big(\frac{(\rho v)^\varepsilon-\rho^\varepsilon v^\varepsilon}{\rho^\varepsilon}\Big) dxdt\B|\\
&\leq C\Big(
\|w^\varepsilon\|_{L^{\f{7}{3(1-\alpha)}}_t(L^{\f{2}{\alpha}}_x)}\|\s v^\varepsilon\|_{L^{\f{2}{\alpha}}_t(L^{\f{4}{2-\alpha}}_x)}
+
\|v^\varepsilon\|_{L^{\f{7}{3(1-\alpha)}}_t(L^{\f{2}{\alpha}}_x)}\|\nabla w^\varepsilon\|_{L^{\f{2}{\alpha}}_t(L^{\f{4}{2-\alpha}}_x)}\\
&~~~
+
\|w^\varepsilon\|_{L^{\f{7}{3(1-\alpha)}}_t(L^{\f{2}{\alpha}}_x)}\|\nabla v^\varepsilon\|_{L^{\f{2}{\alpha}}_t(L^{\f{4}{2-\alpha}}_x)}
\Big)
\left\|
\frac{(\rho v)^\varepsilon-\rho^\varepsilon v^\varepsilon}{\rho^\varepsilon}
\right\|_{L^{\f{14}{8-\alpha}}_t(L^{\f{4}{2-\alpha}}_x)}.
\ea\ee
Due to $v\in L^{\f{2}{\alpha}}(0,T;B^{\alpha}_{\f{4}{2-\alpha},\infty}(\mathbb{T}^3))\cap L^\infty (0,T;L^2(\mathbb{T}^3))$ with $\f23\leq \alpha <1$, we can conclude from the Gagliardo-Nirenberg inequality in the framework of Besov spaces in \cite[Theorem 1.2]{[WWY23]} that 
\be\label{3.17}\|v\|_{L^{\f{7}{3(1-\alpha)}}_t(L^{\f{2}{\alpha}}_x)}\leq C\|v\|_{L^{\f{2}{\alpha}}(0,T;B^{\alpha}_{\f{4}{2-\alpha},\infty})}^{\f{6-6\alpha}{7\alpha}}\|v\|_{L^\infty (0,T;L^2)}^{\f{13\alpha-6}{7\alpha}},\ee
from which it follows that 
\be\label{3.18}
\|\omega^\varepsilon\|_{L^{\f{7}{3(1-\alpha)}}_t(L^{\f{2}{\alpha}}_x)}\leq C\varepsilon^{-1}\|v\|_{L^{\f{7}{3(1-\alpha )}}_t(L^{\f{2}{\alpha}}_x)},
\ee
and 
\be\label{3.19}
\|v^\varepsilon\|_{L^{\f{7}{3(1-\alpha)}}_t(L^{\f{2}{\alpha}}_x)}\leq C\|v\|_{L^{\f{7}{3(1-\alpha )}}_t(L^{\f{2}{\alpha}}_x)}.
\ee
Moreover, applying Lemma \ref{lem2.3} with $\theta=\gamma=1$, $f=v$ and $g=\rho$, we get 
\be\label{3.20}
\ba
\left\|
\frac{(\rho v)^\varepsilon-\rho^\varepsilon v^\varepsilon}{\rho^\varepsilon}
\right\|_{L^{\f{14}{8-\alpha}}_t(L^{\f{4}{2-\alpha}}_x)}\leq &\left\|
(\rho v)^\varepsilon-\rho^\varepsilon v^\varepsilon
\right\|_{L^{\f{14}{8-\alpha}}_t(L^{\f{4}{2-\alpha}}_x)}\\
\leq & o(\varepsilon^{\alpha+\beta})\|\rho\|_{L^{\f{7}{4(1-\alpha)}}_t(B^{\beta}_{\infty,\infty})}\|v\|_{L^{\f{2}{\alpha}}(0,T;B^{\alpha}_{\f{4}{2-\alpha},\infty})}.
\ea\ee
Then gathering \eqref{3.18}-\eqref{3.20} into \eqref{3.16}, it follows that 
\be\label{3.21}\ba
|R_{21}|\leq o(\varepsilon^{\alpha-2+\alpha+\beta}).
\ea\ee
Hence, letting $2\alpha+\beta -2=0$ with $\f{2}{3}<\alpha <1$, which implies  $v\in L^{\f{2}{\alpha}}(0,T;B^{\alpha}_{\f{4}{2-\alpha},c(\mathbb{N})}(\mathbb{T}^3)),\rho \in L^{\f{7}{4(1-\alpha)}}(0,T;B^{2-2\alpha}_{\infty,\infty}(\mathbb{T}^3))$ ensures the $R_{21}$ converges to zero as $\varepsilon \to 0$. 

For term $R_{22}$, taking advantage of Lemma \ref{lem2.3} by letting $\theta=\gamma=1$, $f=v$ and $g=\rho$, we obtain
\be\label{3.22}\ba
|R_{22}|&\leq C\,
\|(\rho v)^\varepsilon-\rho^\varepsilon v^\varepsilon\|_{L^{\f{14}{8-\alpha}}_t(L^{\frac{4}{2-\alpha}}_x)}
\,
\|\nabla v^\varepsilon\|_{L^{\f{7}{3(1-\alpha)}}_t(L^{\f{2}{\alpha}}_x)}
\,
\|
\nabla\Big(\frac{(\rho v)^\varepsilon}{\rho^\varepsilon}\Big)
\|_{L^{\f{2}{\alpha}}_t(L^{\f{4}{2-\alpha}}_x)}\\
&\leq C\varepsilon^{2-2\alpha+\alpha}\|\rho\|_{L^{\f{7}{4(1-\alpha)}}_t(B^{2-2\alpha}_{\infty,\infty})}\|v\|_{L^{\f{2}{\alpha}}(0,T;B^{\alpha}_{\f{4}{2-\alpha},c(\mathbb{N})})}\\
&\quad\varepsilon^{ -1}\| v\|_{L^{\f{7}{3(1-\alpha)}}_t(L^{\f{2}{\alpha}}_x)}\varepsilon^{\alpha -1}\|v\|_{L^{\f{2}{\alpha}}(0,T;B^{\alpha}_{\f{4}{2-\alpha},c(\mathbb{N})})}\\
&\leq o(1).
\ea\ee
Combining \eqref{3.21} and \eqref{3.22}, we can obtain 
\be\label{3.23}
|R_2|\to 0,~\text{as}~\varepsilon \to 0.
\ee
Likewise, proceeding similarly as \eqref{3.23}, we infer that
\be\label{3.24}\ba
|R_3|&\leq C\left\|
(\rho v)^\varepsilon-\rho^\varepsilon v^\varepsilon
\right\|_{L^{\f{14}{8-\alpha}}_t(L^{\frac{4}{2-\alpha}}_x)}\\
&
\left(
\|
\nabla\Big(\frac{(\rho v)^\varepsilon}{\rho^\varepsilon}\Big)
\|_{L^{\f{2}{\alpha}}_t(L^{\f{4}{2-\alpha}}_x)}
+
\|\nabla v^\varepsilon\|_{L^{\f{2}{\alpha}}_t(L^{\f{4}{2-\alpha}}_x)}
\right)\|\partial_t \phi\|_{L^{\f{7}{3(1-\alpha)}}_t(L^{\f{2}{\alpha}}_x)}\\
&\leq o(\varepsilon^{2-2\alpha +\alpha+\alpha-1})\leq o(\varepsilon^1),
\ea\ee
and
\be\label{3.25}\ba
|R_4|&\leq C\,
\|(\rho v)^\varepsilon-\rho^\varepsilon v^\varepsilon\|_{L^{\f{14}{8-\alpha}}_t(L^{\frac{4}{2-\alpha}}_x)}
\,
\|\nabla(\rho v)^\varepsilon\|_{L^{\f{7}{3(1-\alpha)}}_t(L^{\f{2}{\alpha}}_x)}
\,
\|
\nabla\Big(\frac{(\rho v)^\varepsilon}{\rho^\varepsilon}\Big)
\|_{L^{\f{2}{\alpha}}_t(L^{\f{4}{2-\alpha}}_x)}\\
&\leq C\varepsilon^{2-2\alpha+\alpha}\|\rho\|_{L^{\f{7}{4(1-\alpha)}}_t(B^{2-2\alpha}_{\infty,\infty})}\|v\|_{L^{\f{2}{\alpha}}(0,T;B^{\alpha}_{\f{4}{2-\alpha},c(\mathbb{N})})}\\
&\quad\varepsilon^{ -1}\| \rho v\|_{L^{\f{7}{3(1-\alpha)}}_t(L^{\f{2}{\alpha}}_x)}\varepsilon^{\alpha -1}\|v\|_{L^{\f{2}{\alpha}}(0,T;B^{\alpha}_{\f{4}{2-\alpha},c(\mathbb{N})})}\\
&\leq  o(\varepsilon^{2-2\alpha +\alpha-1+\alpha-1})\leq o(1).
\ea\ee

For the pressure term $R_5$, notice that $\pi\in C^2([c_1,c_2])$, it follows that 
\be\label{3.26}\ba
|\pi(s)-\pi(s_0)-\pi^{'}(s_0)(s-s_0)|\leq C|s-s_0|^2.
\ea\ee
Hence, we have
\be\label{3.27}\ba
|\pi(\rho^\varepsilon)-\pi(\rho)-\pi^{'}(\rho)(\rho^\varepsilon-\rho)|\leq C|\rho^\varepsilon-\rho|^2.
\ea\ee
And similarly, we get
\be\label{3.28}\ba
|\pi(\rho(y))-\pi(\rho(x))-\pi^{'}(\rho(x))(\rho(y)-\rho(x))|\leq C|\rho(y)-\rho(x)|^2.
\ea\ee
By virtue of  convolution against the function $\varrho_{\varepsilon}$ with respect to $y$ in \eqref{3.28} and using the Jensen's inequality, it follows that
\be\label{3.29}\ba
|(\pi(\rho))^\varepsilon-\pi(\rho)-\pi^{'}(\rho)(\rho^\varepsilon-\rho)|\leq C|\rho-\rho(\cdot)|^2*\varrho_{\varepsilon}.
\ea\ee
Then in combination with \eqref{3.27} and \eqref{3.29}, we can obtain
\be\label{3.30}\ba
|(\pi(\rho))^\varepsilon-\pi(\rho^\varepsilon)|\leq C|\rho^\varepsilon-\rho|^2+C|\rho -\rho(\cdot)|^2*\varrho_{\varepsilon},
\ea\ee
from which and  $\rho\in L^{\f{7}{4(1-\alpha)}}_t(B^{\beta}_{\infty,\infty})$ it follows that
\begin{align*}
	\|(\pi(\rho))^\varepsilon-\pi(\rho^\varepsilon)\|_{L^{\f{7}{8(1-\alpha)}}_t(L^\infty_x)}
	\le
	C\varepsilon^{2\beta}
	\|\rho\|^2_{L^{\f{7}{4(1-\alpha)}}_t(B^{\beta}_{\infty,\infty})}.
\end{align*}
As a consequence, integrating by parts in $R_5$ and using $\s\big(\nabla\times\frac{(\rho v)^\varepsilon}{\rho^\varepsilon}\big)=0$ together with H\"older's inequality, we obtain
\be\label{3.31}\ba
|R_5|&\leq 
C
\|(\pi(\rho))^\varepsilon-\pi(\rho^\varepsilon)\|_{L^{\f{7}{8(1-\alpha)}}_t(L^\infty_x)}
\left\|\nabla\Big(\tfrac1{\rho^\varepsilon}\Big)\right\|_{L^{\f{7}{4(1-\alpha)}}_t(L^\infty_x)}\\
&~~
\left\|\nabla\Big(\tfrac{(\rho v)^\varepsilon}{\rho^\varepsilon}\Big)\right\|_{L^{\f{2}{\alpha}}_t(L^{\f{4}{2-\alpha}}_x)}\|\phi(t)\|_{L^{\f{14}{17\alpha -10}}_t(L^{\f{4}{\alpha+2}}_x)}\\
&\le \varepsilon^{2\beta}
\|\rho\|^2_{L^{\f{7}{4(1-\alpha)}}_t(B^{\beta}_{\infty,\infty})}\varepsilon^{\beta -1}
\|\rho\|_{L^{\f{7}{4(1-\alpha)}}_t(B^{\beta}_{\infty,\infty})} \varepsilon^{\alpha -1}\|v\|_{L^{\f{2}{\alpha}}(0,T;B^{\alpha}_{\f{4}{2-\alpha},c(\mathbb{N})})}\\
&\leq o(\varepsilon ^{3\beta+\alpha-2})\to 0,
\ea\ee
provided $\beta =\f{2-\alpha}{3}$. Hence, together with the previous constraint for $\rho$, we need $\rho \in L^{\f{7}{4(1-\alpha)}}(0,T; B^{\max\{2-2\alpha, \f{2-\alpha}{3}\}}_{\infty,\infty}(\mathbb{T}^3))\cap L^\infty (0,T;L^\infty(\mathbb{T}^3))$. 
Then, collecting \eqref{3.13},\eqref{3.23}, \eqref{3.24}, \eqref{3.25} and \eqref{3.31}, the Lebesgue dominated convergence theorem yields that the helicity is locally conserved in the sense of distributions of time, namely,
\be\label{3.32}\ba
&\f{1}{2}\int_0^T\int_{\mathbb T^3} v \cdot \omega \partial_t \phi(t) \, dxdt=0.
\ea\ee

{\bf Endpoint case: $\alpha=\f{2}{3}$.} 

Now, we are in a position to prove the terms $R_1-R_5$ converge to zero under the constraint that 
\be\label{3.33}\ba
&0<c_1\leq \rho \leq c_2<\infty,~\rho \in L^{\f{21}{4}}(0,T;B^{\f{2}{3}}_{\infty,\infty}(\mathbb{T}^3)),\\
&\rho v,v \in L^{3}(0,T;B^{\f{2}{3}}_{3,c(\mathbb{N})}(\mathbb{T}^3)).
\ea\ee
Firstly,  by H\"older's inequality, we arrive at
\be\label{3.34}\begin{aligned} |R_1|&=\B|\int_0^T\int_{\mathbb{T}^3}\phi(t)
\Big[(\rho v\otimes v)^\varepsilon-(\rho v)^\varepsilon\otimes v^\varepsilon\Big]:
\B[\!\left(\nabla\times(\frac{(\rho v)^\varepsilon}{\rho^\varepsilon})\right)\otimes
\nabla\!(\frac{1}{\rho^\varepsilon})\\
&\quad
+
\frac{1}{\rho^\varepsilon}\,
\nabla\!\left(\nabla\times(\frac{(\rho v)^\varepsilon}{\rho^\varepsilon})\right)\B]\,dxdt\B|\\
&\le
C\,
\|(\rho v\otimes v)^\varepsilon-(\rho v)^\varepsilon\otimes v^\varepsilon\|_{L^{\frac{3}{2}}_t(L^{\frac{3}{2}}_x)}  \\
&\quad
\Big(
\|(\nabla\times(\frac{(\rho v)^\varepsilon}{\rho^\varepsilon}))\otimes\nabla(\tfrac1{\rho^\varepsilon})\|_{L^3_t(L^3_x)}
+
\|\tfrac1{\rho^\varepsilon}\nabla(\nabla\times\tfrac{(\rho v)^\varepsilon}{\rho^\varepsilon})\|_{L^3_t(L^3_x)}
\Big).
\\&\le
C\,
\|(\rho v\otimes v)^\varepsilon-(\rho v)^\varepsilon\otimes v^\varepsilon\|_{L^{\frac{3}{2}}_t(L^{\frac{3}{2}}_x)} \\
&\quad
\Big(
\|\nabla(\tfrac1{\rho^\varepsilon})\|_{L^\infty_t(L^\infty_x)}
\left\|\nabla\Big(\tfrac{(\rho v)^\varepsilon}{\rho^\varepsilon}\Big)\right\|_{L^3_t(L^3_x)}
+
\left\|\nabla^2\Big(\tfrac{(\rho v)^\varepsilon}{\rho^\varepsilon}\Big)\right\|_{L^3_t(L^3_x)}
\Big)\\
&\leq C\varepsilon^{\f{4}{3}}\|\rho v\|_{L^3_t({B}^{\frac{2}{3}}_{3,\infty})}\|v\|_{L^3_t(B^{\frac{2}{3}}_{3,c(\mathbb{N})})}\B(\varepsilon^{ -1}\|\rho\|_{L^\infty_t(L^\infty_x)}\varepsilon^{\frac{2}{3} -1}\|v\|_{L^3_t(B^{\frac{2}{3}}_{3,c(\mathbb{N})})}\\
&\quad+\varepsilon^{\frac{2}{3} -2}\|v\|_{L^3_t(B^{\frac{2}{3}}_{3,c(\mathbb{N})})}\B)\\
&\leq o(1).
\ea\ee
Next, Recall \eqref{3.14}
\be\label{3.35}\ba
|R_2|
&\le 
\B|\int_0^T \int_{\mathbb{T}^3}\phi(t)
\Big(\nabla\times\frac{(\rho v)^\varepsilon}{\rho^\varepsilon}\Big)\cdot
\big(\nabla v^\varepsilon\cdot v^\varepsilon\big)
\,dxdt\B|
\\
&\quad+\B|
\int_0^T \int_{\mathbb{T}^3} \phi(t)
\Big(\nabla\times\frac{(\rho v)^\varepsilon}{\rho^\varepsilon}\Big)\cdot
\Big(\nabla v^\varepsilon\cdot \frac{(\rho v)^\varepsilon-\rho^\varepsilon v^\varepsilon}{\rho^\varepsilon}\Big)
\,dxdt\B|\\
&=R_{21}+R_{22}.
\end{aligned}\ee
Before proceeding further, the Gagliardo-Nirenberg inequality in  Besov spaces in \eqref{3.17} implies that 
\be\label{3.36}\|v\|_{L^{7}_t(L^{3}_x)}\leq C\|v\|_{L^{3}(0,T;B^{\f{2}{3}}_{3,\infty})}^{\f{3}{7}}\|v\|_{L^\infty (0,T;L^2)}^{\f{4}{7}}.\ee
Then, combining the   vector identity
$
\nabla\times(w^\varepsilon\times v^\varepsilon)=w^\varepsilon\s v^\varepsilon - v^\varepsilon\s w^\varepsilon + v^\varepsilon\cdot\nabla w^\varepsilon - w^\varepsilon\cdot\nabla v^\varepsilon,
$
and  $\s w^\varepsilon=0$, the first term $R_{21}$ can be controlled as 
\be\label{3.37}\ba
|R_{21}|
&= \B|\int_0^T \int_{\mathbb{T}^3}\phi(t)
\nabla\times(\omega^\varepsilon\times v^\varepsilon)
\cdot
\Big(\frac{(\rho v)^\varepsilon-\rho^\varepsilon v^\varepsilon}{\rho^\varepsilon}\Big)
\,dxdt\B|\\
&\leq \B|\int_0^T\int_{\mathbb{T}^3} \phi(t)\left(w^\varepsilon\s v^\varepsilon  + v^\varepsilon\cdot\nabla w^\varepsilon - w^\varepsilon\cdot\nabla v^\varepsilon\right)\cdot \Big(\frac{(\rho v)^\varepsilon-\rho^\varepsilon v^\varepsilon}{\rho^\varepsilon}\Big) dxdt\B|\\
&\leq C\Big(
\|w^\varepsilon\|_{L^7_t(L^3_x)}\|\s v^\varepsilon\|_{L^3_t(L^3_x)}
+
\|v^\varepsilon\|_{L^7_t(L^3_x)}\|\nabla w^\varepsilon\|_{L^3_t(L^3_x)}
+
\|w^\varepsilon\|_{L^7_t(L^3_x)}\|\nabla v^\varepsilon\|_{L^3_t(L^3_x)}
\Big)\\
&\quad\times\left\|
\frac{(\rho v)^\varepsilon-\rho^\varepsilon v^\varepsilon}{\rho^\varepsilon}
\right\|_{L^{\frac{21}{11}}_t(L^3_x)}
\\
&\leq C \left(\| \nabla v^\varepsilon \|_{L^7_t(L^3_x)} \| \nabla v^\varepsilon \|_{L^3_t(L^3_x)} + \| v^\varepsilon \|_{L^7_t(L^3_x)} \| \nabla^2 v^\varepsilon \|_{L^3_t(L^3_x)} \right)\\
&\quad\times\| (\rho v)^\varepsilon - \rho^\varepsilon v^\varepsilon \|_{L^{\frac{21}{11}}_t(L^3_x)}\\
&\leq C\left(\varepsilon^{\frac{2}{3}-2}\| v \|_{L^7_t(L^3_x)}\|v\|_{L^3_t(B^{\frac{2}{3}}_{3,c(\mathbb{N})})}+\varepsilon^{\frac{2}{3} -2}\| v \|_{L^7_t(L^3_x)}\|v\|_{L^3_t(B^{\frac{2}{3}}_{3,c(\mathbb{N})})}\right)\\
&\quad\times\varepsilon^{\frac{4}{3}}\|\rho\|_{L^{\f{21}{4}}_t(B^{\frac{2}{3}}_{\infty,\infty})}\|v\|_{L^3_t(B^{\frac{2}{3}}_{3,c(\mathbb{N})})}\\
&\leq o(1).
\ea\ee
By virtue of H\"older's inequality, Lemma \ref{lem2.5} and Lemma \ref{lem2.4}, we arrive at
\be\label{3.38}\ba
|R_{22}|&\leq C\,
\|(\rho v)^\varepsilon-\rho^\varepsilon v^\varepsilon\|_{L^3_t(L^3_x)}
\,
\|\nabla v^\varepsilon\|_{L^3_t(L^3_x)}
\,
\|
\nabla\Big(\frac{(\rho v)^\varepsilon}{\rho^\varepsilon}\Big)
\|_{L^3_t(L^3_x)}\\
&\leq C\varepsilon^{\frac{2}{3}}\|\rho\|_{L^\infty_t(L^\infty_x)}\|v\|_{L^3_t(B^{\frac{2}{3}}_{3,c(\mathbb{N})})}\varepsilon^{\frac{2}{3} -1}\| v\|_{L^3_t(B^{\frac{2}{3}}_{3,c(\mathbb N)})}\varepsilon^{\frac{2}{3} -1}\|v\|_{L^3_t(B^{\frac{2}{3}}_{3,c(\mathbb{N})})}\\
&\leq o(1).
\ea\ee
Making use of  H\"older's inequality, Lemma \ref{lem2.1}, Lemma \ref{lem2.5} and Lemma \ref{lem2.4}, we infer that
\be\label{3.39}\ba
|R_3|&\leq C\left\|
(\rho v)^\varepsilon-\rho^\varepsilon v^\varepsilon
\right\|_{L^3_t(L^3_x)}
\left(
\|
\nabla\Big(\frac{(\rho v)^\varepsilon}{\rho^\varepsilon}\Big)
\|_{L^3_t(L^3_x)}
+
\|\nabla v^\varepsilon\|_{L^3_t(L^3_x)}
\right)\|\partial_t \phi(t)\|_{L^3_t(L^3_x)}\\
&\leq C\varepsilon^{\frac{2}{3}}\|\rho\|_{L^\infty_t(L^\infty_x)}\|v\|_{L^3_t(B^{\frac{2}{3}}_{3,c(\mathbb{N})})}\varepsilon^{\frac{2}{3} -1}\|v\|_{L^3_t(B^{\frac{2}{3}}_{3,c(\mathbb{N})})}\\
&\leq o(\varepsilon^{\frac{1}{3}}),
\ea\ee
and
\be\label{3.40}\ba
|R_4|&\leq C\,
\|(\rho v)^\varepsilon-\rho^\varepsilon v^\varepsilon\|_{L^3_t(L^3_x)}
\,
\|\nabla(\rho v)^\varepsilon\|_{L^3_t(L^3_x)}
\,
\|
\nabla\Big(\frac{(\rho v)^\varepsilon}{\rho^\varepsilon}\Big)
\|_{L^3_t(L^3_x)}\\
&\leq C\varepsilon^{\frac{2}{3}}\|\rho\|_{L^\infty_t(L^\infty_x)}\|v\|_{L^3_t(B^{\frac{2}{3}}_{3,c(\mathbb{N})})}\varepsilon^{\frac{2}{3} -1}\|\rho v\|_{L^3_t(B^{\frac{2}{3}}_{3,\infty})}\varepsilon^{\frac{2}{3}-1}\|v\|_{L^3_t(B^{\frac{2}{3}}_{3,c(\mathbb{N})})}\\
&\leq o(1).
\ea\ee
Furthermore,
for term $R_5$, integrating by parts as in \eqref{3.31}, it follows from the H\"older's inequality that
\be\label{3.41}\ba
|R_5|
&\le
C\,
\|(\pi(\rho))^\varepsilon-\pi(\rho^\varepsilon)\|_{L^{\f{21}{4}}_t(L^\infty_x)}
\left\|\nabla\Big(\tfrac1{\rho^\varepsilon}\Big)\right\|_{L^{\f{21}{4}}_t(L^\infty_x)}
\left\|\nabla\Big(\tfrac{(\rho v)^\varepsilon}{\rho^\varepsilon}\Big)\right\|_{L^3_t(L^3_x)}\|\phi(t)\|_{L^{\f{7}{2}}_t(L^{\f32}_x)}.
\ea\ee
Since $\pi\in C^1([c_1,c_2])$, the mean value theorem yields that
\be\label{3.42}\ba
\pi(s)-\pi(s_0)&=\int_{s_0}^{s} \pi^{'}(\tau) d\tau=\int_{s_0}^s \pi^{'}({s_0}) d\tau+\int_{s_0}^s \pi^{'}(\tau)-\pi^{'}(s_0) d\tau\\
&=\pi^{'}(s_0)(s-s_0)+\int_{s_0}^s \pi^{'}(\tau)-\pi^{'}(s_0) d\tau,
\ea\ee
which implies
\be\label{3.43}\ba
|\pi(s)-\pi(s_0)-\pi^{'}(s_0)(s-s_0)|\leq C|s-s_0|.
\ea\ee
Hence, we have
\be\label{3.44}\ba
|\pi(\rho^\varepsilon)-\pi(\rho)-\pi^{'}(\rho)(\rho^\varepsilon-\rho)|\leq C|\rho^\varepsilon-\rho|.
\ea\ee
And similarly, we get
\be\label{3.45}\ba
|\pi(\rho(y))-\pi(\rho(x))-\pi^{'}(\rho(x))(\rho(y)-\rho(x))|\leq C|\rho(y)-\rho(x)|.
\ea\ee
By virtue of  convolution against the function $\varrho_{\varepsilon}$ with respect to $y$ in \eqref{3.45} and using the Jensen's inequality, it follows that
\be\label{3.46}\ba
|(\pi(\rho))^\varepsilon-\pi(\rho)-\pi^{'}(\rho)(\rho^\varepsilon-\rho)|\leq C|\rho-\rho(\cdot)|*\varrho_{\varepsilon}.
\ea\ee
Then in combination with \eqref{3.44} and \eqref{3.46}, we can obtain
\be\label{3.47}\ba
|(\pi(\rho))^\varepsilon-\pi(\rho^\varepsilon)|\leq C|\rho^\varepsilon-\rho|+C|\rho -\rho(\cdot)|*\varrho_{\varepsilon},
\ea\ee
from which and  $\rho \in L^{\f{21}{4}}(0,T;B^{\f{2}{3}}_{\infty,\infty}(\mathbb{T}^3))$ it follows that
\begin{align*}
	\|(\pi(\rho))^\varepsilon-\pi(\rho^\varepsilon)\|_{L^{\f{21}{4}}_t(L^\infty_x)}
	\le
	C\varepsilon^{\frac{2}{3}}
	\|\rho\|_{L^{\f{21}{4}}_t(B^{\frac{2}{3}}_{\infty,\infty})},
\end{align*}
Consequently, we have
\be\label{3.48}\ba
|R_5|
&\le \varepsilon^{\frac{2}{3}}
\|\rho\|_{L^{\f{21}{4}}_t(B^{\frac{2}{3}}_{\infty,\infty})}\varepsilon^{\frac{2}{3} -1}
\|\rho\|_{L^{\f{21}{4}}_t(B^{\frac{2}{3}}_{\infty,\infty})} \varepsilon^{\frac{2}{3} -1}\|v\|_{L^3_t(B^{\frac{2}{3}}_{3,c(\mathbb{N})})}\\
&\leq o(1).
\ea\ee
Then, collecting \eqref{3.34}-\eqref{3.35}, \eqref{3.37}-\eqref{3.40} and \eqref{3.48}, then the Lebesgue dominated convergence theorem yields that the helicity is locally conserved in the sense of distributions of time, namely,
\be\label{3.49}\ba
&\f{1}{2}\int_0^T\int_{\mathbb T^3} v \cdot \omega \partial_t \phi(t) \, dxdt=0.
\ea\ee

{\bf Condition (iii) of Theorem~\ref{the 1.1}.}
 
 Now, we prove that the terms $R_1$--$R_5$ converge to zero under the following assumptions
\be\label{3.50}
\ba
&0<c_1\leq \rho \leq c_2<\infty,\quad \rho \in L^{\infty}(0,T;B^{\frac{1}{3}}_{\infty,c(\mathbb{N})}(\mathbb{T}^3)),\\
&v\in L^{2}(0,T;W^{1,4}(\mathbb{T}^3)),\quad \pi\in C^2([c_1,c_2]).
\ea
\ee
Using Young's convolution inequality, $\|\nabla\varrho_\varepsilon\|_{L^1(\mathbb{R}^3)}\leq C\varepsilon^{-1}$, and \eqref{3.50} , we obtain
\be\label{3.51}\ba
&\|\nabla v^\varepsilon\|_{L^\infty_t(L^2_x)}+\|\omega^\varepsilon\|_{L^\infty_t(L^2_x)}+\|\nabla(\rho v)^\varepsilon\|_{L^\infty_t(L^2_x)}\leq C\varepsilon^{-1},\\
&\|\nabla\omega^\varepsilon\|_{L^2_t(L^4_x)}\leq C\varepsilon^{-1}\|\nabla v\|_{L^2_t(L^4_x)}.
\ea\ee
Applying Lemma~\ref{lem2.1}  to $f=\rho$ with $\alpha=\frac{1}{3}$, $p=q=\infty$ and $k=1$, and using $\rho^\varepsilon\geq c_1$, we obtain
\be\label{3.52}
\Big\|\nabla\Big(\frac{1}{\rho^\varepsilon}\Big)\Big\|_{L^\infty_t(L^\infty_x)}\leq c_1^{-2}\|\nabla\rho^\varepsilon\|_{L^\infty_t(L^\infty_x)}\leq C\varepsilon^{-\frac{2}{3}}.
\ee
Moreover, applying Lemma~\ref{lem-sob} (i) and (ii) with $p=2$, $q=4$ and $\beta=\frac{1}{3}$, we know that
\be\label{3.53}
\Big\|\nabla\Big(\frac{(\rho v)^\varepsilon}{\rho^\varepsilon}\Big)\Big\|_{L^2_t(L^4_x)}\leq C\|\nabla v\|_{L^2_t(L^4_x)},\qquad
\lim_{\varepsilon\to0}\varepsilon\Big\|\nabla^2\Big(\frac{(\rho v)^\varepsilon}{\rho^\varepsilon}\Big)\Big\|_{L^2_t(L^4_x)}=0.
\ee
Furthermore, using the Sobolev translation estimate and the triangle inequality, we have, for a.e. $t\in(0,T)$,
\begin{equation*}\ba
&\frac{\|v(\cdot-y)-v(\cdot)\|_{L^4}}{|y|^{\frac{5}{6}}}\leq |y|^{\frac{1}{6}}\|\nabla v\|_{L^4},\quad 0<|y|\leq1,\\
&\frac{\|v(\cdot-y)-v(\cdot)\|_{L^4}}{|y|^{\frac{5}{6}}}\leq\frac{\|v(\cdot-y)\|_{L^4}+\|v\|_{L^4}}{|y|^{\frac{5}{6}}}
=\frac{2\|v\|_{L^4}}{|y|^{\frac{5}{6}}}\leq2\|v\|_{L^4},\quad |y|>1.
\ea\end{equation*}
Using \eqref{2.2}, \eqref{2.3} and
$v\in L^2(0,T;W^{1,4}(\mathbb{T}^3))$, we obtain
\begin{equation*}\ba
	\|v\|_{L^2(0,T;B^{\frac{5}{6}}_{4,c(\mathbb{N})}(\mathbb{T}^3))}
	&=\left(\int_0^T\left[\|v(t)\|_{L^4}
	+\sup_{|y|>0}\frac{\|v(t,\cdot-y)-v(t,\cdot)\|_{L^4}}
	{|y|^{\frac{5}{6}}}\right]^2dt\right)^{\frac{1}{2}}\\
	&\leq\left(\int_0^T
	\left[3\|v(t)\|_{L^4}+\|\nabla v(t)\|_{L^4}\right]^2dt
	\right)^{\frac{1}{2}}\\
	&\leq C\|v\|_{L^2(0,T;W^{1,4}(\mathbb{T}^3))}<\infty.
	\ea\end{equation*}
Here, the $c(\mathbb{N})$ condition in \eqref{2.3} follows from,
for a.e. $t\in(0,T)$ and $0<h\leq1$,
\begin{equation*}
	\sup_{0<|y|\leq h}
	\frac{\|v(t,\cdot-y)-v(t,\cdot)\|_{L^4}}{|y|^{\frac{5}{6}}}
	\leq h^{\frac{1}{6}}\|\nabla v(t)\|_{L^4}
	\to0\quad(h\to0).
\end{equation*}
 Applying Lemma~\ref{lem2.3} to $f=v$ and $g=\rho$ with $\theta=\gamma=1$, $\alpha=\frac{5}{6}$, $\beta=\frac{1}{3}$, $(p_1,q_1)=(2,4)$ and $(p_3,q_3)=(\infty,\infty)$, we obtain
\be\label{3.54}
\|(\rho v)^\varepsilon-\rho^\varepsilon v^\varepsilon\|_{L^2_t(L^4_x)}\leq o(\varepsilon^{\frac{7}{6}}).
\ee

Firstly, taking $p=2$ and $q=4$ in \eqref{3.9}, we obtain
\be\label{3.55}\ba
|R_1|\leq{}& C\|(\rho v\otimes v)^\varepsilon-(\rho v)^\varepsilon\otimes v^\varepsilon\|_{L^2_t(L^{\frac{4}{3}}_x)}\\
&\times\Big(\Big\|\nabla\Big(\frac{(\rho v)^\varepsilon}{\rho^\varepsilon}\Big)\Big\|_{L^2_t(L^4_x)}\Big\|\nabla\Big(\frac{1}{\rho^\varepsilon}\Big)\Big\|_{L^\infty_t(L^\infty_x)}+c_1^{-1}\Big\|\nabla^2\Big(\frac{(\rho v)^\varepsilon}{\rho^\varepsilon}\Big)\Big\|_{L^2_t(L^4_x)}\Big).
\ea\ee
Then applying \cite[Lemma 2.2]{[WY]} componentwise with $f=v_j\in L^{2}(0,T;W^{1,4}(\mathbb{T}^3))$, $g=\rho v_i\in L^{\infty}(0,T;L^{2}(\mathbb{T}^3))$ and $(p_1,q_1,p_2,q_2)=(2,4,\infty,2)$, we obtain
\be\label{3.56}
\|(\rho v\otimes v)^\varepsilon-(\rho v)^\varepsilon\otimes v^\varepsilon\|_{L^2_t(L^{\frac{4}{3}}_x)}\leq C\varepsilon\|v\|_{L^2_t(W^{1,4}_x)}\|\rho v\|_{L^\infty_t(L^2_x)}\leq C\varepsilon.
\ee
Then, combining \eqref{3.55} and \eqref{3.56} with \eqref{3.52} and \eqref{3.53}, we obtain
\be\label{3.57}
|R_1|\leq C\varepsilon^{\frac{1}{3}}\|\nabla v\|_{L^2_t(L^4_x)}+C\varepsilon\Big\|\nabla^2\Big(\frac{(\rho v)^\varepsilon}{\rho^\varepsilon}\Big)\Big\|_{L^2_t(L^4_x)}\to0\quad(\varepsilon\to0).
\ee

Next, recalling \eqref{3.14}, we have $|R_2|\leq R_{21}+R_{22}$. Applying H\"older's inequality to the integral for $R_{21}$ in \eqref{3.16}, and using  \eqref{3.51} and \eqref{3.54}, we obtain
\be\label{3.58}\ba
|R_{21}|&\leq C\Big(\|\omega^\varepsilon\|_{L^\infty_t(L^2_x)}\|\nabla v^\varepsilon\|_{L^2_t(L^4_x)}+\|v^\varepsilon\|_{L^\infty_t(L^2_x)}\|\nabla\omega^\varepsilon\|_{L^2_t(L^4_x)}\Big)\\
&\quad\times\|(\rho v)^\varepsilon-\rho^\varepsilon v^\varepsilon\|_{L^2_t(L^4_x)}\\
&\leq C\varepsilon^{-1}\|\nabla v\|_{L^2_t(L^4_x)}\,o(\varepsilon^{\frac{7}{6}})\to0\quad(\varepsilon\to0).
\ea\ee
For the term $R_{22}$ in \eqref{3.14}, applying H\"older's inequality and using \eqref{3.51}, \eqref{3.53} and \eqref{3.54}, we obtain
\be\label{3.59}\ba
|R_{22}|&\leq c_1^{-1}\Big\|\nabla\times\frac{(\rho v)^\varepsilon}{\rho^\varepsilon}\Big\|_{L^2_t(L^4_x)}\|\nabla v^\varepsilon\|_{L^\infty_t(L^2_x)}\|(\rho v)^\varepsilon-\rho^\varepsilon v^\varepsilon\|_{L^2_t(L^4_x)}\\
&\leq C\varepsilon^{-1}\|\nabla v\|_{L^2_t(L^4_x)}\,o(\varepsilon^{\frac{7}{6}})\to0\quad(\varepsilon\to0).
\ea\ee

Likewise, applying H\"older's inequality to $R_3$ and $R_4$, and using \eqref{3.51}, \eqref{3.53} and \eqref{3.54}, we obtain
\be\label{3.60}\ba
|R_3|&\leq C\|(\rho v)^\varepsilon-\rho^\varepsilon v^\varepsilon\|_{L^2_t(L^4_x)}\Big(\|\omega^\varepsilon\|_{L^2_t(L^4_x)}+\Big\|\nabla\times\frac{(\rho v)^\varepsilon}{\rho^\varepsilon}\Big\|_{L^2_t(L^4_x)}\Big)\|\partial_t\phi\|_{L^\infty_t(L^2_x)}\\
&\leq C\|\nabla v\|_{L^2_t(L^4_x)}\,o(\varepsilon^{\frac{7}{6}})\to0,\\
|R_4|&\leq C\|(\rho v)^\varepsilon-\rho^\varepsilon v^\varepsilon\|_{L^2_t(L^4_x)}\|\s(\rho v)^\varepsilon\|_{L^\infty_t(L^2_x)}\Big\|\nabla\times\frac{(\rho v)^\varepsilon}{\rho^\varepsilon}\Big\|_{L^2_t(L^4_x)}\\
&\leq C\varepsilon^{-1}\|\nabla v\|_{L^2_t(L^4_x)}\,o(\varepsilon^{\frac{7}{6}})\to0\quad(\varepsilon\to0).
\ea\ee

For the term $R_5$, integrating by parts in $R_5$ and using $\s\big(\nabla\times\frac{(\rho v)^\varepsilon}{\rho^\varepsilon}\big)=0$, we obtain
\be\label{3.61}
R_5=-\int_0^T\int_{\mathbb{T}^3}\phi(t)\big((\pi(\rho))^\varepsilon-\pi(\rho^\varepsilon)\big)\nabla\Big(\frac{1}{\rho^\varepsilon}\Big)\cdot\nabla\times\Big(\frac{(\rho v)^\varepsilon}{\rho^\varepsilon}\Big)\,dxdt.
\ee
Since $\pi\in C^2([c_1,c_2])$, using \eqref{3.30} and $|\rho^\varepsilon-\rho|\leq\int_{B_\varepsilon(0)}\varrho_\varepsilon(y)|\rho(x-y)-\rho(x)|\,dy$, we obtain, for a.e. $t\in(0,T)$,
\be\label{3.62}\ba
\|(\pi(\rho))^\varepsilon-\pi(\rho^\varepsilon)\|_{L^\infty_x}&\leq C\sup_{0<|y|<\varepsilon}\|\rho(t,\cdot-y)-\rho(t,\cdot)\|^2_{L^\infty}\\
&\leq C\varepsilon^{\frac{2}{3}}\Big(\sup_{0<|y|<\varepsilon}\frac{\|\rho(t,\cdot-y)-\rho(t,\cdot)\|_{L^\infty}}{|y|^{\frac{1}{3}}}\Big)^2.
\ea\ee
Using \eqref{2.3} and \eqref{3.50}, we have, for a.e. $t\in(0,T)$,
\begin{equation*}\ba
&\sup_{0<|y|<\varepsilon}\frac{\|\rho(t,\cdot-y)-\rho(t,\cdot)\|_{L^\infty}}{|y|^{\frac{1}{3}}}
\leq C\|\rho(t)\|_{B^{\frac{1}{3}}_{\infty,\infty}},\\
&\sup_{0<|y|<\varepsilon}\frac{\|\rho(t,\cdot-y)-\rho(t,\cdot)\|_{L^\infty}}{|y|^{\frac{1}{3}}}
\to0\quad(\varepsilon\to0).
\ea\end{equation*}
Hence, applying the Lebesgue dominated convergence theorem to \eqref{3.62}, we obtain
\be\label{3.63}\ba
\|(\pi(\rho))^\varepsilon-\pi(\rho^\varepsilon)\|_{L^2_t(L^\infty_x)}
&\leq C\varepsilon^{\frac{2}{3}}\left(\int_0^T\left(\sup_{0<|y|<\varepsilon}
\frac{\|\rho(t,\cdot-y)-\rho(t,\cdot)\|_{L^\infty}}{|y|^{\frac{1}{3}}}\right)^4\,dt\right)^{\frac{1}{2}}\\
&=o(\varepsilon^{\frac{2}{3}})\quad(\varepsilon\to0).
\ea\ee
Consequently, applying H\"older's inequality to \eqref{3.61} and using \eqref{3.52}, \eqref{3.53} and \eqref{3.63}, we obtain
\be\label{3.64}\ba
|R_5|&\leq C\|(\pi(\rho))^\varepsilon-\pi(\rho^\varepsilon)\|_{L^2_t(L^\infty_x)}\Big\|\nabla\Big(\frac{1}{\rho^\varepsilon}\Big)\Big\|_{L^\infty_t(L^\infty_x)}\Big\|\nabla\times\frac{(\rho v)^\varepsilon}{\rho^\varepsilon}\Big\|_{L^2_t(L^4_x)}\\
&\leq o(\varepsilon^{\frac{2}{3}})\,C\varepsilon^{-\frac{2}{3}}\|\nabla v\|_{L^2_t(L^4_x)}\to0\quad(\varepsilon\to0).
\ea\ee
Then, collecting \eqref{3.57}--\eqref{3.60} and \mbox{\eqref{3.64}},  we let $\varepsilon\to0$ in \eqref{3.8} to obtain
\begin{equation*}
\frac{1}{2}\int_0^T\int_{\mathbb T^3} v \cdot \omega\, \partial_t \phi(t) \, dxdt=0.
\end{equation*}

{
{\bf Condition (iv) of Theorem~\ref{the 1.1}.} In this case, we aim to relax the constraint on the density $\rho $
 by imposing an additional constraint on the momentum $\rho v$.
 
We now assume that
\be\label{3.65}
\ba
&0<c_1\leq \rho \leq c_2<\infty,\quad \rho \in L^{\infty}(0,T;B^{\frac{1}{3}}_{\infty,\infty}(\mathbb{T}^3)),\\
&v\in L^{2}(0,T;W^{1,4}(\mathbb{T}^3)),~\rho v\in L^1(0,T;W^{1,1}(\mathbb{T}^3)),~\pi\in C^2(0,\infty).
\ea
\ee
We remark that the estimates \eqref{3.51}--\eqref{3.60} remain valid under assumption  \eqref{3.65}, as  their derivations rely only on $\rho\in L^\infty_t(B^{\frac{1}{3}}_{\infty,\infty})$ and $v\in L^2_t(W^{1,4}_x)$ in addition to the common hypotheses of Theorem~\ref{the 1.1}. Thus, it suffices to show that  $R_5\to 0,$ as $\varepsilon\to0$.

Recall \eqref{3.61} and take advantage of the vector identity $E\cdot (E\times F)=(E\times E)\cdot F=0$ , we have
\be\label{3.66}\ba
R_5&=-\int_0^T\int_{\mathbb{T}^3}\phi(t)\big((\pi(\rho))^\varepsilon-\pi(\rho^\varepsilon)\big)\nabla\Big(\frac{1}{\rho^\varepsilon}\Big)\cdot\Big[\frac{\big(\nabla\times(\rho v)\big)^\varepsilon}{\rho^\varepsilon}+\nabla\Big(\frac{1}{\rho^\varepsilon}\Big)\times(\rho v)^\varepsilon\Big]\,dxdt\\
&=-\int_0^T\int_{\mathbb{T}^3}\phi(t)\frac{(\pi(\rho))^\varepsilon-\pi(\rho^\varepsilon)}{\rho^\varepsilon}\,\nabla\Big(\frac{1}{\rho^\varepsilon}\Big)\cdot\big(\nabla\times(\rho v)\big)^\varepsilon\,dxdt.
\ea\ee
Using \eqref{3.62}, \eqref{2.2} and \eqref{3.65}, we infer that, for a.e. $t\in(0,T)$,
\be\label{3.67}\ba
\|(\pi(\rho))^\varepsilon-\pi(\rho^\varepsilon)\|_{L^\infty_x}&\leq C\sup_{0<|y|<\varepsilon}\|\rho(t,\cdot-y)-\rho(t,\cdot)\|^2_{L^\infty}\leq C\varepsilon^{\frac{2}{3}}\|\rho\|^2_{L^\infty_t(B^{\frac{1}{3}}_{\infty,\infty})},
\ea\ee
which together with \eqref{3.52} and \eqref{3.67} yields  
\be\label{3.68}\ba
\Big\|\frac{(\pi(\rho))^\varepsilon-\pi(\rho^\varepsilon)}{\rho^\varepsilon}\,\nabla\Big(\frac{1}{\rho^\varepsilon}\Big)\Big\|_{L^\infty_t(L^\infty_x)}&\leq c_1^{-1}\|(\pi(\rho))^\varepsilon-\pi(\rho^\varepsilon)\|_{L^\infty_t(L^\infty_x)}\Big\|\nabla\Big(\frac{1}{\rho^\varepsilon}\Big)\Big\|_{L^\infty_t(L^\infty_x)}\\
&\leq C\varepsilon^{\frac{2}{3}}\varepsilon^{-\frac{2}{3}}=C.
\ea\ee
As a consequence, we have
\be\label{3.68-1}
\ba
|R_5|
&\leq \left| -\int_0^T\int_{\mathbb{T}^3}\phi(t)\frac{(\pi(\rho))^\varepsilon-\pi(\rho^\varepsilon)}{\rho^\varepsilon}\,\nabla\Big(\frac{1}{\rho^\varepsilon}\Big)\cdot\big(\nabla\times(\rho v)\big)^\varepsilon\,dxdt\right|\\
&\leq\|\phi\|_{L^\infty(0,T)}
\Big\|\frac{(\pi(\rho))^\varepsilon-\pi(\rho^\varepsilon)}{\rho^\varepsilon}\,
\nabla\Big(\frac{1}{\rho^\varepsilon}\Big)\Big\|_{L^\infty_t(L^\infty_x)}
\|\nabla(\rho v)\|_{L^1(L^1)}\\
&\leq C.
\ea\ee
On the other hand, by virtue of  \eqref{3.66}, \eqref{3.68} and H\"older's inequality, we arrive at
\be\label{3.69}\ba
|R_5|&=\Big|\int_0^T\int_{\mathbb{T}^3}\phi(t)\frac{(\pi(\rho))^\varepsilon-\pi(\rho^\varepsilon)}{\rho^\varepsilon}\,\nabla\Big(\frac{1}{\rho^\varepsilon}\Big)\cdot\Big[\nabla\times(\rho v)+\big(\nabla\times(\rho v)\big)^\varepsilon-\nabla\times(\rho v)\Big]\,dxdt\Big|\\
&\leq\Big|\int_0^T\int_{\mathbb{T}^3}\phi(t)\frac{(\pi(\rho))^\varepsilon-\pi(\rho^\varepsilon)}{\rho^\varepsilon}\,\nabla\Big(\frac{1}{\rho^\varepsilon}\Big)\cdot\nabla\times(\rho v)\,dxdt\Big|\\
&\quad+\Big|\int_0^T\int_{\mathbb{T}^3}\phi(t)\frac{(\pi(\rho))^\varepsilon-\pi(\rho^\varepsilon)}{\rho^\varepsilon}\,\nabla\Big(\frac{1}{\rho^\varepsilon}\Big)\cdot\Big[\big(\nabla\times(\rho v)\big)^\varepsilon-\nabla\times(\rho v)\Big]\,dxdt\Big|\\
&\leq\Big|\int_0^T\int_{\mathbb{T}^3}\phi(t)\frac{(\pi(\rho))^\varepsilon-\pi(\rho^\varepsilon)}{\rho^\varepsilon}\,\nabla\Big(\frac{1}{\rho^\varepsilon}\Big)\cdot\nabla\times(\rho v)\,dxdt\Big|\\
&\quad+C\Big\|\frac{(\pi(\rho))^\varepsilon-\pi(\rho^\varepsilon)}{\rho^\varepsilon}\,\nabla\Big(\frac{1}{\rho^\varepsilon}\Big)\Big\|_{L^\infty_t(L^\infty_x)}\big\|\big(\nabla\times(\rho v)\big)^\varepsilon-\nabla\times(\rho v)\big\|_{L^1_t(L^1_x)}\\
&=\Big|\int_0^T\int_{\mathbb{T}^3}\phi(t)\frac{(\pi(\rho))^\varepsilon-\pi(\rho^\varepsilon)}{\rho^\varepsilon}\,\nabla\Big(\frac{1}{\rho^\varepsilon}\Big)\cdot\nabla\times(\rho v)\,dxdt\Big|+o(1)\quad(\varepsilon\to0).
\ea\ee
\begin{samepage}
Here, applying H\"older's inequality and using \eqref{3.65}, we obtain
\begin{equation*}\ba
\big\|\big(\nabla\times(\rho v)\big)^\varepsilon-\nabla\times(\rho v)\big\|_{L^1_t(L^1_x)}\leq C
\big\|\big(\nabla(\rho v)\big)^\varepsilon-\nabla(\rho v)\big\|_{L^1_t(L^1_x)}\to0
\quad(\varepsilon\to0),
\ea\end{equation*}
which justifies the $o(1)$ term in \eqref{3.69}.
\end{samepage}

We next estimate the first term on the right-hand side of \eqref{3.69}. Applying Sobolev's embedding theorem and using \eqref{3.65}, we obtain, for a.e. $t\in(0,T)$,
\be\label{3.70}
v(t,\cdot)\in W^{1,4}(\mathbb{T}^3)\hookrightarrow C(\mathbb{T}^3).
\ee
For fixed time $t$, we let $U_t:=\{x\in\mathbb{T}^3:v(t,x)\neq0\}$. Combining $\rho\leq c_2$ and \eqref{3.70}, we  conclude 
\begin{equation*}
x\in\mathbb{T}^3\setminus U_t
\quad\Longrightarrow\quad v(t,x)=0
\quad\Longrightarrow\quad |(\rho v)(t,x)|\leq c_2|v(t,x)|=0.
\end{equation*}
Applying Lemma~\ref{lem-zero-set} componentwise to $(\rho v)(t,\cdot)$, we obtain
\be\label{3.71}
\nabla(\rho v)(t,\cdot)=0,\qquad
\nabla\times(\rho v)(t,\cdot)=0
\quad\text{a.e. on }\mathbb{T}^3\setminus U_t.
\ee
Thus, for $x_0\in U_t$, choose $j\in\{1,2,3\}$ such that $v_j(t,x_0)\neq0$. By \eqref{3.70}, there exists $r>0$ such that
\be\label{3.72}
0<c_0:=\frac{1}{2}|v_j(t,x_0)|\leq |v_j(t,x)|,
\qquad x\in B_{2r}(x_0)\subset U_t.
\ee
Then the Sobolev product and chain rules together with \eqref{3.72}  implies 
\begin{equation*}\ba
\rho(t,\cdot)&=\frac{(\rho v_j)(t,\cdot)}{v_j(t,\cdot)}
\in W^{1,4}(B_{2r}(x_0)),\\
\nabla\rho(t,\cdot)&=\nabla\left(\frac{(\rho v_j)(t,\cdot)}{v_j(t,\cdot)}\right)=\frac{v_j(t,\cdot)\nabla(\rho v_j)(t,\cdot)-(\rho v_j)(t,\cdot)\nabla v_j(t,\cdot)}{v_j(t,\cdot)^2}\\&\qquad\qquad\qquad\qquad=\frac{\nabla(\rho v_j)(t,\cdot)-\rho(t,\cdot)\nabla v_j(t,\cdot)}{v_j(t,\cdot)},
\ea\end{equation*}
which further yields
\begin{equation*}
\|\nabla\rho(t)\|_{L^1(B_{2r}(x_0))}
\leq c_0^{-1}\Big(\|\nabla(\rho v_j)(t)\|_{L^1(\mathbb{T}^3)}
+c_2\|\nabla v_j(t)\|_{L^1(\mathbb{T}^3)}\Big)<\infty.
\end{equation*}
Thus we conclude that $\rho(t,\cdot)\in W^{1,1}_{\mathrm{loc}}(U_t)$.  Moreover, 
Since $B_\varepsilon(x)\subset B_{2r}(x_0)$ for $x\in B_r(x_0)$ and $0<\varepsilon<r$, integrating by parts, we obtain, at every Lebesgue point $x\in B_r(x_0)$ of $\nabla\rho(t,\cdot)$,
\begin{equation*}\ba
|\nabla\rho^\varepsilon(t,x)-\nabla\rho(t,x)|
&=\Big|\int_{B_\varepsilon(0)}\varrho_\varepsilon(y)
\big(\nabla\rho(t,x-y)-\nabla\rho(t,x)\big)\,dy\Big|\\
&\leq\int_{B_\varepsilon(0)}\varrho_\varepsilon(y)
|\nabla\rho(t,x-y)-\nabla\rho(t,x)|\,dy\\
&\leq\frac{C}{|B_\varepsilon(x)|}\int_{B_\varepsilon(x)}
|\nabla\rho(t,y)-\nabla\rho(t,x)|\,dy\to0\quad(\varepsilon\to0).
\ea\end{equation*}
Then the Lebesgue differentiation theorem yields
\begin{equation*}
\nabla\rho^\varepsilon(t,x)\to\nabla\rho(t,x)
\quad\text{for a.e. }x\in B_r(x_0),~ \text{as}~~ \varepsilon\to0.
\end{equation*}
 In light of \eqref{3.70}, $U_t$ is open, and hence it is covered by countably many such balls $B_r(x_0)$, which gives $\nabla\rho^\varepsilon(t,x)\to\nabla\rho(t,x)$ for a.e. $x\in U_t$ as $\varepsilon\to0$.

Hence, we obtain that, for a.e. $x\in U_t$ and all sufficiently small $\varepsilon$,
\begin{equation*}
|\nabla\rho^\varepsilon(t,x)-\nabla\rho(t,x)|\leq1.
\end{equation*}
Using \eqref{3.67}, $\rho^\varepsilon\geq c_1$ and
$\nabla\big(\frac{1}{\rho^\varepsilon}\big)=-\frac{\nabla\rho^\varepsilon}{(\rho^\varepsilon)^2}$, we obtain
\be\label{3.73}\ba
\Big|\frac{(\pi(\rho))^\varepsilon(t,x)-\pi(\rho^\varepsilon(t,x))}{\rho^\varepsilon(t,x)}\,
\nabla\Big(\frac{1}{\rho^\varepsilon}\Big)(t,x)\Big|&=\frac{|(\pi(\rho))^\varepsilon(t,x)-\pi(\rho^\varepsilon(t,x))|}{(\rho^\varepsilon(t,x))^3}|\nabla\rho^\varepsilon(t,x)|\\
&\leq Cc_1^{-3}\varepsilon^{\frac{2}{3}}|\nabla\rho^\varepsilon(t,x)|\\
&=Cc_1^{-3}\varepsilon^{\frac{2}{3}}
\Big|\nabla\rho(t,x)
+\big(\nabla\rho^\varepsilon(t,x)-\nabla\rho(t,x)\big)\Big|\\
&\leq C\varepsilon^{\frac{2}{3}}
\Big(|\nabla\rho(t,x)|
+|\nabla\rho^\varepsilon(t,x)-\nabla\rho(t,x)|\Big)\\
&\leq C\varepsilon^{\frac{2}{3}}\big(|\nabla\rho(t,x)|+1\big)\to0
\\
&\qquad\text{for a.e. }x\in U_t\quad(\varepsilon\to0).
\ea\ee
Combining \eqref{3.71} and \eqref{3.73} and applying Fubini's theorem, we obtain
\begin{equation*}\ba
&\phi(t)\frac{(\pi(\rho))^\varepsilon-\pi(\rho^\varepsilon)}{\rho^\varepsilon}\,
\nabla\Big(\frac{1}{\rho^\varepsilon}\Big)\cdot\nabla\times(\rho v)\to0, \quad \text{a.e. in }(0,T)\times\mathbb{T}^3\quad(\varepsilon\to0).
\ea\end{equation*}
Thus in combination with \eqref{3.68-1} and \eqref{3.73}, applying the Lebesgue dominated convergence theorem, we obtain
\be\label{3.74}
\Big|\int_0^T\int_{\mathbb{T}^3}\phi(t)\frac{(\pi(\rho))^\varepsilon-\pi(\rho^\varepsilon)}{\rho^\varepsilon}\,\nabla\Big(\frac{1}{\rho^\varepsilon}\Big)\cdot\nabla\times(\rho v)\,dxdt\Big|\to0
\quad(\varepsilon\to0).
\ee
Combining \eqref{3.69} and \eqref{3.74}, we conclude that $R_5\to0$ as $\varepsilon\to0$.

}

Then, collecting \eqref{3.57}--\eqref{3.60} and {the corresponding limit $R_5\to0$ in case (iii) or (iv)}, and using $v^\varepsilon\to v$ in $C([0,T];L^{3}(\mathbb{T}^3))$ and $\omega^\varepsilon\to\omega$ in $C([0,T];L^{\frac{3}{2}}(\mathbb{T}^3))$, we let $\varepsilon\to0$ in \eqref{3.8} to obtain
\be\label{3.75}
\frac{1}{2}\int_0^T\int_{\mathbb T^3} v \cdot \omega\, \partial_t \phi(t) \, dxdt=0.
\ee

\underline{Step 2:} Next, we  aim to get the helicity equality up to the initial time $t=0$. Before we go any further, it should be noted that \eqref{3.32} (as well as \eqref{3.49} and \eqref{3.75}) remains valid for function $\phi$ belonging to $W^{1,\infty}$ rather than $C^1$, then for any $t_0\in(0,T)$, we redefine the test function $\phi$ as $\phi_\tau$ for some positive $\tau$ and $\eta$ such that $\tau +\eta <t_0$ and $t_0+\eta<T$, that is
\begin{equation*}
	\phi_\tau(t)=\left\{\begin{array}{lll}
		0, & 0\leq t\leq \tau,\\
		\f{t-\tau}{\eta}, & \tau\leq t\leq \tau+\eta,\\
		1, &\tau+\eta\leq t\leq t_0,\\
		\f{t_0+\eta-t}{\eta}, & t_0\leq t\leq t_0 +\eta,\\
		0, & t_0+\eta\leq t.
	\end{array}\right.
\end{equation*}
Then substituting this test function into \eqref{3.32}, we arrive at
\begin{equation*}
	\f{1}{2}\int_\tau^{\tau+\eta}\int_{\mathbb{T}^3} \f{1}{\eta}\left(v\cdot\omega \right)dxdt-\f{1}{2}\int_{t_0}^{t_0+\eta}\int_{\mathbb{T}^3} \f{1}{\eta}\left(v\cdot\omega \right)dxdt=0.
\end{equation*}
Taking $\eta\rightarrow 0$ and using  the continuity of $t\mapsto\int_{\mathbb{T}^3}v\cdot\omega\, dx$ on $[0,T]$, which follows from $v\in C([0,T];L^{3}(\mathbb{T}^3))$ and $\omega\in C([0,T];L^{\f{3}{2}}(\mathbb{T}^3))$, we deduce that
\begin{equation}
	\begin{aligned}
		\f{1}{2}\int_{\mathbb{T}^3} \left(v\cdot\omega \right)(t_0) dx=	\f{1}{2}\int_{\mathbb{T}^3} \left(v\cdot\omega \right)(\tau) dx.
	\end{aligned}
\end{equation}
Finally, letting $\tau\rightarrow 0$ and using the continuity again, we can obtain
\begin{equation}\ba
		\f{1}{2}\int_{\mathbb{T}^3} \left(v\cdot\omega \right)(t_0) dx=	\f{1}{2}\int_{\mathbb{T}^3} \left(v\cdot\omega \right)(0) dx.
	\ea\end{equation}
Then we have completed the proof of Theorem \ref{the 1.1}.

\end{proof}

\section*{Acknowledgements}
Ye was
partially sponsored by the Training Plan for Young Backbone Teachers in Colleges and Universities of Henan Province (No. 2024GGJS022).  \\

{\bf Conflict of Interest Statement}

This work does not have any conflicts of interest.

{\bf Data Availability}

This publication is supported by multiple datasets, which are openly available at locations cited in the reference section.

\end{document}